\documentclass[twocolumn, reqno]{amsart}
\usepackage{xcolor}
\usepackage{amsfonts,amssymb}
\usepackage[margin=18mm,,footskip=7mm]{geometry}
\usepackage[dvipsnames]{xcolor}   
\usepackage{hyperref} % clickable URL
\hypersetup{
    colorlinks=true,
    linkcolor=OliveGreen,  
    urlcolor=BlueViolet,
    citecolor=MidnightBlue
    }
\usepackage{tabularray}
\usepackage{microtype}

\DefTblrTemplate{contfoot-text}{default}{continued}
\DefTblrTemplate{conthead}{default}{conthead}
\DefTblrTemplate{capcont}{default}{}

\newtheorem{thm}{Theorem}

\begin{document}

\title[Simple groups with narrow spectrum]{Simple groups with narrow prime spectrum: Extended list}

\author{Andrei V. Zavarnitsine}%
\address{Andrei V. Zavarnitsine
\newline\indent Sobolev Institute of Mathematics,
\newline\indent 4, Koptyug av.
\newline\indent 630090, Novosibirsk, Russia
} 
\email{zav@math.nsc.ru
\newline\indent ORCID: \href{https:\\orcid.org/0000-0003-1983-3304}{0000-0003-1983-3304}
}
\thanks{This research was carried out within the State Contract of the Sobolev Institute of Mathematics (FWNF-2026-0017).}
\maketitle

\renewenvironment{quote}
  {\list{}{\leftmargin=9mm \rightmargin=9mm}%
   \item\relax}
  {\endlist}

\begin{quote}
\noindent{\sc Abstract.} Generalising a previous result, we determine 
all non-abelian finite simple groups whose order has largest prime divisor not exceeding $10^4$. The computer code for this and similar  calculations is made available.


\medskip
\noindent{\small {\sc MSC2020:} 
20D60 % Arithmetic and combinatorial problems
}\hfill {\small {\sc UDC:}  512.542.6}

\end{quote}

\bigskip
\section*{Introduction}
This paper continues the work of~\cite{Z09} which has proved useful. 
In~\cite{Z09}, all non-abelian finite simple groups with order having largest prime divisor not exceeding~$1000$ were determined. Here we extend this bound to~$10000$, thereby proving the following. 

\begin{thm}\label{main} There are $15072$ isomorphism types of finite non-abelian simple groups whose order has all prime divisors less than $10000$.
\end{thm}

Recall from~\cite{Z09} that,
for a finite group~$G$, its \emph{prime spectrum}, denoted by $\pi(G)$, is the set of prime divisors of~$|G|$.
Given a prime~$p$, we write $\mathfrak{S}_p$ for the set of non-abelian finite simple groups~$G$ satisfying 
\[
p \in \pi(G) \subseteq \{2,3,5,\ldots,p\}.
\]
The sets $\mathfrak{S}_p$, $p \geqslant 5$, are finite and always contain the \emph{generic} groups
\[
  L_2(p),\; A_p,\; A_{p+1},\; \ldots,\; A_{p'-1},
\]
where $p'$ is the smallest prime exceeding~$p$. In particular, there are  \(p' - p + 1\) 
such groups in $\mathfrak{S}_p$.
The primes~$p$ for which $\mathfrak{S}_p$ consists solely of generic groups are called
\emph{generic primes}. The non-generic elements of~$\mathfrak{S}_p$, when they exist, are of particular interest.

The $13100$ groups obtained in the present paper (which, together with the $1972$ already listed in~\cite{Z09}, give a total of $15072$ groups) can be inferred from Tables~\ref{gen} and~\ref{ng}. In Table~\ref{gen}, we list the $301$ generic primes between $1000$ and $10000$. In Table~\ref{ng}, we list the $3041$ non-generic groups~$G$ 
from the union 
\[
  \mathfrak{S}_{1009}\cup\ldots\cup \mathfrak{S}_{9973}
\]
of the sets $\mathfrak{S}_p$ corresponding to the $760$ remaining non-generic primes between $1000$ and $10000$. 

The notation for simple groups in Table~\ref{ng} follows that of \cite{Atlas}, except the sizes of defining fields bigger than $100$ are written exponentially,
e.\,g., $L_2(3^7)$ in place of $L_2(2187)$.

As is apparent from the tables, the overwhelming majority of groups are generic. It is reasonable to expect alternating groups to dominate the set $\mathfrak{S}_p$ asymptotically.

Curiously, the largest number $44$ of non-generic groups in~$\mathfrak{S}_p$ for $p<1000$ is attained at $p = 257$ (a Fermat prime) and remains so for all $p<10000$.
We have no example of~$\mathfrak{S}_p$ with more than~$44$ non-generic elements. The second largest value, $39$, occurs at $p = 6481$ which is the greatest prime factor of~$3^{12}+1$.

The calculations of this paper (and, retrospectively, those of~\cite{Z09}) can be verified using the \textsf{GAP}~\cite{GAP} code provided in~\cite{ZG}. Specifically, we implement a universal function 
\begin{flalign}\label{sgpi}
  & \texttt{SimpleGroupsPi(pi)} &
\end{flalign}
which accepts an arbitrary set of primes \texttt{pi} and returns (codes of) all non-abelian finite simple groups~$G$ satisfying $\pi(G) \subseteq \texttt{pi}$. 

The idea behind the implemented algorithm is due to V.\,D.~Mazurov and is outlined in \cite[p.~51]{m94}. Roughly, any simple group $G$ with $\pi(G) \subseteq \pi$  satisfies $G\in\mathcal{S}\cup\mathcal{A}_\pi\cup\mathcal{L}_\pi$, where $\mathcal{S}$ is the set of $26$
sporadic groups, $\mathcal{A}_\pi$ is the set of alternating groups $A_n$ with $n\leqslant p_0-1$, where
$p_0$ is the smallest prime greater than all primes in $\pi$, and $\mathcal{L}_\pi$ is the set of
groups of Lie type defined over a field of order $p^k$ with $p\in \pi$
and having rank~$l_p$, where  
$$k\leqslant t_p=\max_{r\in \pi\setminus\{p\}}\operatorname{ord}_r p,$$
and $l_p\leqslant \max\{8,t_p\}$. Hence, finding all such groups $G$ reduces to checking the orders of finitely many groups.

The running time of function (\ref{sgpi}) depends on both the cardinality of \texttt{pi} and the size of its largest element: the computations with all primes up to~$1000$ take approximately $1$~minute, while those with primes less than~$10000$ require $25$~hours. It is therefore impractical to apply this linear implementation directly to, say, all primes up to~$10^5$.
However, the algorithm is highly parallelisable.

On the other hand, the program runs quickly when \texttt{pi} is small, even if it contains large primes. Thus, in \cite[Example~1]{ZG}, we find all $13$ simple groups 
\begin{multline*}
A_{5},\ A_{6},\ U_{5}(2),\ L_{2}(3^{5}),\ S_{4}(3),\ L_{2}(11),\ L_{2}(11^{2}),\\
S_{4}(11),\ U_{3}(11),\ U_{4}(11),\ U_{5}(11),\ M_{11},\ M_{12}
\end{multline*}
with prime spectrum a subset of
\[
\pi(U_5(11))=\{2,\ 3,\ 5,\ 11,\ 37,\ 61,\ 13421\}
\]
within seconds, where $U_5(11)$ is notable for having the largest prime divisor of its order, $13421$, among all groups listed in the Atlas table~\cite[p.~239]{Atlas}. Calculations of this sort with a custom set \texttt{pi} might be useful for applications.

Theorem \ref{main} results from running (\ref{sgpi}) with 
\[
 \texttt{pi}=\{2,3,\ldots,9973\}.
\]
Tables~\ref{gen} and~\ref{ng} are obtained after sorting the output according to the maximal prime factor and determining whether \(\mathfrak{S}_p\) is generic, see \cite[Example 3]{ZG} for details.

Another natural stratification of the groups from Theorem \ref{main} is
by the size of their prime spectrum. 
Excluding the alternating groups, the largest size of prime spectrum turns out to be $24$ attained for the four groups 
\begin{equation}\label{four}
L_{15}(4),\ S_{30}(2),\ O_{32}^+(2),\ U_{15}(4),
\end{equation}
of which the first three have the same prime spectrum equal to
\begin{equation}\label{pi24}
\begin{split}
  \{ 2, 3, 5, 7, 11,& 13, 17, 19, 23, 29, 31, 41, 43, 73, 89,\\ 
  113,& 127, 151, 241, 257, 331, 683, 2731, 8191 \}.
\end{split}
\end{equation}
The second largest size, $23$, is attained for just two non-alternating groups, $O_{30}^+(2)$ and $O_{30}^-(2)$,
whose prime spectrum can be obtained from (\ref{pi24}) by removing $331$ and $151$, respectively.

More generally, for every $n=3,\ldots,24$, we list in \cite[Example 4]{ZG} explicitly the sets $\mathfrak{K}_n=\mathfrak{K}_{n,1229}$ of all {\em $n$-primary} groups (also known as {\em $K_n$-groups}) from Theorem~\ref{main}. Here $\mathfrak{K}_{n,m}$ denotes the set of non-abelian simple groups whose orders have exactly $n$ distinct prime factors all of which do not exceed the $m$-th prime $p_m$. We have $\mathfrak{K}_{n,m_1}\subseteq \mathfrak{K}_{n,m_2}$ whenever $m_1\leqslant m_2$. Note that $1229$ is the index of $9973$, the largest prime less than $10000$.

The sizes of $\mathfrak{K}_n$ for $n\leqslant 24$ are as follows:

\begin{gather*}
\begin{array}{cccccccc}
n                 & 3 &  4 &   5 &   6 &   7 &   8 &   9    \\ \hline
|\mathfrak{K}_n|^{\vphantom{A^A}}  & 8 & 65 & 349 & 715 & 595 & 628 & 828  
\end{array} \\
\begin{array}{cccccccc}
n                 &  10 &  11 &  12  & 13  &  14 & 15 & 16    \\ \hline
|\mathfrak{K}_n|^{\vphantom{A^A}}  & 641 & 398 & 309 & 238  & 119 & 87 & 62   
\end{array} \\
\begin{array}{ccccccccc}
n                 & 17 & 18 & 19 & 20 & 21 & 22 & 23 & 24   \\ \hline
|\mathfrak{K}_n|^{\vphantom{A^A}}  & 32 & 12 & 18 &  7 & 16 & 15 &  8 & 12  
\end{array}
\end{gather*}

\medskip
\noindent
Clearly, $\mathfrak{K}_n$ contains the $p_{n+1}-p_n$ alternating groups 
\[
A_{p_n}, \ A_{p_n+1}, \ldots,  A_{p_{n+1}-1}
\]
for every $n=3,\ldots,1229$ and, in fact, only those groups for $n>24$. For example, \(\mathfrak{K}_{24}\) consists of $A_{89},\ldots,A_{96}$ plus the four groups in (\ref{four}), a total of $12$ groups.

\medskip
{\em Acknowledgement.} The author is expressly grateful to Prof.\ D.\,O. Revin for helpful discussions of this topic.

\newpage
\section*{The tables}

\NewTblrTheme{narrow}{
\SetTblrStyle{foot}{\itshape}
\SetTblrStyle{caption-tag}{font=\bfseries}
}

\begin{longtblr}[
  theme = narrow,
  caption = {Primes $p\in\{1000,\ldots,10000\}$ with generic $\mathfrak{S}_p$},
  label = {gen}
]{
  colspec = {c},
}
\hline
1009, 1013, 1019, 1033, 1039, 1097, 1103, 1151, 1163, \\
1187, 1193, 1213, 1217, 1249, 1259, 1279, 1307, 1319, \\
1361, 1381, 1409, 1439, 1453, 1481, 1523, 1559, 1579, \\
1627, 1667, 1669, 1721, 1733, 1777, 1811, 1847, 1879, \\
1901, 1907, 1933, 1949, 1997, 2003, 2011, 2029, 2063, \\
2069, 2087, 2129, 2137, 2179, 2221, 2239, 2341, 2351, \\
2357, 2377, 2381, 2399, 2423, 2447, 2459, 2477, 2543, \\
2549, 2593, 2647, 2659, 2663, 2699, 2711, 2741, 2857, \\
2879, 2887, 2909, 2927, 2963, 3023, 3061, 3067, 3119, \\
3163, 3167, 3191, 3209, 3217, 3271, 3299, 3329, 3371, \\
3407, 3469, 3491, 3539, 3677, 3697, 3719, 3733, 3761, \\
3767, 3797, 3803, 3821, 3823, 3847, 3877, 3911, 3923, \\
3929, 3943, 3989, 4001, 4019, 4127, 4139, 4153, 4157, \\
4159, 4231, 4241, 4259, 4339, 4349, 4397, 4409, 4447, \\
4451, 4463, 4507, 4517, 4547, 4583, 4663, 4703, 4759, \\
4871, 4919, 4931, 4943, 4993, 4999, 5003, 5011, 5021, \\
5039, 5081, 5147, 5179, 5273, 5279, 5297, 5303, 5309, \\
5333, 5387, 5393, 5399, 5407, 5417, 5449, 5471, 5501, \\
5519, 5521, 5563, 5651, 5669, 5689, 5737, 5741, 5783, \\
5839, 5897, 5903, 5939, 5953, 5981, 5987, 6011, 6047, \\
6101, 6113, 6131, 6199, 6203, 6221, 6271, 6329, 6361, \\
6563, 6599, 6653, 6689, 6709, 6719, 6761, 6791, 6857, \\
6869, 6883, 6947, 7001, 7013, 7019, 7043, 7069, 7109, \\
7211, 7229, 7243, 7247, 7331, 7349, 7393, 7411, 7417, \\
7451, 7457, 7487, 7517, 7523, 7541, 7547, 7573, 7589, \\
7591, 7643, 7691, 7741, 7757, 7793, 7817, 7823, 7829, \\
7901, 7937, 7951, 8017, 8053, 8059, 8081, 8087, 8123, \\
8167, 8171, 8231, 8233, 8293, 8311, 8353, 8419, 8431, \\
8467, 8513, 8543, 8627, 8677, 8693, 8699, 8741, 8753, \\
8803, 8819, 8849, 8963, 8969, 9013, 9029, 9041, 9059, \\
9137, 9151, 9203, 9293, 9311, 9319, 9323, 9371, 9397, \\
9413, 9437, 9473, 9497, 9521, 9533, 9539, 9547, 9613, \\
9623, 9629, 9719, 9743, 9749, 9767, 9769, 9781, 9787, \\
9803, 9829, 9857, 9973\\
\hline
\end{longtblr}

\newpage

\begin{longtblr}[
  theme = narrow,
  caption = {Non-generic simple groups $G$ with \\  $p\in\pi(G)\subseteq\{2,3,5,\ldots,p\}$ for $1000<p<10000$},
  label = {ng}
]{colspec = {l|l|Q[l,co=1]},
  rowhead = 1
}
 \hline
\SetRow{c} $p$ & $|\mathfrak{S}_p|$ &  $G$\\
\hline
1021 & 16 & $L_{2}(647^{2})$, $S_{4}(647)$, $L_{2}(653^{3})$, $G_2(653)$, $U_{3}(653)$\\
\SetRow{lightgray!40} 1031 & 4 & $U_{3}(1031)$\\
1049 & 4 & $L_{3}(1049)$\\
\SetRow{lightgray!40} 1051 & 15 & $L_{2}(181^{3})$, $G_2(181)$, $U_{3}(181)$, $U_{3}(1051)$\\
1061 & 6 & $L_{2}(103^{2})$, $S_{4}(103)$, $U_{4}(103)$\\
\SetRow{lightgray!40} 1063 & 12 & $L_{3}(7^{3})$, $L_{4}(7^{3})$, $L_{3}(719)$, $L_{2}(719^{3})$, $G_2(719)$\\
1069 & 21 & $U_{3}(983)$, $L_{3}(1069)$\\
\SetRow{lightgray!40} 1087 & 28 & $L_{3}(257)$, $L_{2}(257^{3})$, $G_2(257)$, $L_{3}(829)$, $L_{4}(829)$, $L_{3}(829^{2})$, $L_{2}(829^{3})$, $S_{6}(829)$, $O_{7}(829)$, $O^+_{8}(829)$, $G_2(829)$, $L_{3}(1087)$, $L_{4}(1087)$, $L_{2}(1087^{2})$, $L_{3}(1087^{2})$, $L_{2}(1087^{3})$, $S_{4}(1087)$, $S_{6}(1087)$, $O_{7}(1087)$, $O^+_{8}(1087)$, $G_2(1087)$, $U_{3}(1087)$, $U_{4}(1087)$\\
1091 & 4 & $U_{3}(1091)$\\
\SetRow{lightgray!40} 1093 & 34 & $L_{7}(3)$, $L_{8}(3)$, $L_{9}(3)$, $L_{10}(3)$, $L_{7}(9)$, $L_{2}(3^{7})$, $L_{8}(9)$, $L_{9}(9)$, $S_{14}(3)$, $S_{16}(3)$, $S_{18}(3)$, $O_{15}(3)$, $O_{17}(3)$, $O_{19}(3)$, $O^+_{14}(3)$, $O^+_{16}(3)$, $O^+_{18}(3)$, $O^+_{20}(3)$, $O^-_{16}(3)$, $O^-_{18}(3)$, $E_7(3)$, $L_{3}(151)$, $L_{4}(151)$, $L_{2}(563^{2})$, $S_{4}(563)$, $U_{4}(563)$, $L_{3}(941)$, $L_{2}(941^{3})$, $G_2(941)$\\
1109 & 21 & $L_{3}(1109)$, $L_{4}(1109)$, $L_{2}(1109^{2})$, $L_{3}(1109^{2})$, $L_{2}(1109^{3})$, $S_{4}(1109)$, $S_{6}(1109)$, $O_{7}(1109)$, $O^+_{8}(1109)$, $G_2(1109)$, $U_{3}(1109)$, $U_{4}(1109)$\\
\SetRow{lightgray!40} 1117 & 16 & $L_{2}(11^{6})$, $S_{4}(11^{3})$, $G_2(11^{2})$, $^3D_4(11)$, $U_{3}(11^{2})$, $U_{3}(23^{3})$, $L_{2}(997^{3})$, $G_2(997)$, $U_{3}(997)$\\
1123 & 9 & $L_{2}(1123^{2})$, $S_{4}(1123)$\\
\SetRow{lightgray!40} 1129 & 33 & $L_{4}(31^{2})$, $L_{2}(31^{4})$, $S_{8}(31)$, $S_{4}(31^{2})$, $O_{9}(31)$, $O^-_{8}(31)$, $L_{3}(1129)$, $L_{2}(1129^{3})$, $G_2(1129)$, $U_{3}(1129)$\\
1153 & 15 & $U_{3}(503)$, $L_{2}(1013^{2})$, $S_{4}(1013)$, $L_{3}(1153)$\\
\SetRow{lightgray!40} 1171 & 14 & $U_{3}(421)$, $U_{4}(421)$, $U_{3}(751)$\\
1181 & 22 & $L_{10}(9)$, $L_{5}(3^{4})$, $L_{2}(3^{10})$, $S_{20}(3)$, $S_{10}(9)$, $S_{4}(3^{5})$, $O_{21}(3)$, $O_{11}(9)$, $O^+_{12}(9)$, $O^-_{10}(9)$, $O^-_{20}(3)$, $O^-_{22}(3)$, $U_{5}(9)$, $U_{6}(9)$, $U_{4}(3^{5})$\\
\SetRow{lightgray!40} 1201 & 37 & $L_{4}(49)$, $L_{2}(7^{4})$, $L_{3}(7^{4})$, $L_{2}(7^{12})$, $S_{8}(7)$, $S_{4}(49)$, $S_{6}(49)$, $S_{4}(7^{6})$, $O_{9}(7)$, $O_{7}(49)$, $O^+_{8}(49)$, $O^-_{8}(7)$, $O^-_{10}(7)$, $G_2(7^{4})$, $F_4(7)$, $U_{8}(7)$, $U_{4}(49)$, $U_{3}(7^{4})$, $L_{2}(571^{3})$, $G_2(571)$, $U_{3}(571)$, $L_{2}(631^{3})$, $G_2(631)$, $U_{3}(631)$\\
1223 & 9 & $L_{2}(1223^{2})$, $S_{4}(1223)$\\
\SetRow{lightgray!40} 1229 & 7 & $L_{2}(1229^{2})$, $S_{4}(1229)$, $U_{3}(1229)$, $U_{4}(1229)$\\
1231 & 9 & $U_{3}(127)$, $L_{3}(1231)$\\
\SetRow{lightgray!40} 1237 & 17 & $L_{2}(691^{2})$, $S_{4}(691)$, $U_{3}(937)$, $L_{3}(1237)$\\
1277 & 8 & $L_{4}(113)$, $L_{2}(113^{2})$, $S_{4}(113)$, $L_{2}(1277^{2})$, $S_{4}(1277)$\\
\SetRow{lightgray!40} 1283 & 8 & $L_{3}(1283)$\\
1289 & 5 & $L_{2}(479^{2})$, $S_{4}(479)$\\
\SetRow{lightgray!40} 1291 & 9 & $U_{3}(347)$, $L_{3}(1291)$\\
1297 & 9 & $L_{2}(1297^{2})$, $S_{4}(1297)$, $U_{3}(1297)$, $U_{4}(1297)$\\
\SetRow{lightgray!40} 1301 & 5 & $L_{2}(1301^{2})$, $S_{4}(1301)$\\
1303 & 7 & $L_{2}(1303^{2})$, $S_{4}(1303)$\\
\SetRow{lightgray!40} 1321 & 28 & $L_{5}(2^{12})$, $L_{2}(2^{30})$, $S_{10}(64)$, $S_{4}(2^{15})$, $O^+_{12}(64)$, $O^-_{10}(64)$, $G_2(2^{10})$, $^3D_4(32)$, $U_{5}(64)$, $U_{3}(2^{10})$, $U_{6}(64)$, $L_{4}(257)$, $L_{2}(257^{2})$, $L_{3}(257^{2})$, $S_{4}(257)$, $S_{6}(257)$, $O_{7}(257)$, $O^+_{8}(257)$, $U_{4}(257)$, $Sz(2^{15})$, $^2F_4(32)$\\
1327 & 42 & $L_{3}(347)$, $L_{2}(347^{3})$, $G_2(347)$, $L_{2}(1327^{2})$, $S_{4}(1327)$, $U_{3}(1327)$, $U_{4}(1327)$\\
\SetRow{lightgray!40} 1367 & 8 & $L_{3}(1367)$\\
1373 & 11 & $L_{2}(1373^{2})$, $S_{4}(1373)$\\
\SetRow{lightgray!40} 1399 & 13 & $U_{3}(1009)$, $U_{3}(1399)$\\
1423 & 6 & $L_{3}(643)$\\
\SetRow{lightgray!40} 1427 & 5 & $L_{2}(1427^{2})$, $S_{4}(1427)$\\
1429 & 9 & $L_{4}(809)$, $L_{2}(809^{2})$, $S_{4}(809)$, $L_{3}(1429)$\\
\SetRow{lightgray!40} 1433 & 11 & $L_{2}(1433^{2})$, $S_{4}(1433)$, $U_{3}(1433)$, $U_{4}(1433)$\\
1447 & 7 & $U_{3}(743)$, $U_{4}(743)$\\
\SetRow{lightgray!40} 1451 & 4 & $L_{3}(1451)$\\
1459 & 14 & $L_{3}(1459)$\\
\SetRow{lightgray!40} 1471 & 14 & $L_{3}(251)$, $L_{4}(251)$, $L_{3}(1471)$\\
1483 & 7 & $L_{2}(1483^{2})$, $S_{4}(1483)$\\
\SetRow{lightgray!40} 1487 & 7 & $L_{2}(1487^{2})$, $S_{4}(1487)$, $U_{3}(1487)$, $U_{4}(1487)$\\
1489 & 6 & $U_{3}(1489)$\\
\SetRow{lightgray!40} 1493 & 13 & $L_{2}(1061^{2})$, $S_{4}(1061)$, $L_{3}(1493)$, $L_{4}(1493)$, $L_{2}(1493^{2})$, $S_{4}(1493)$\\
1499 & 14 & $U_{3}(1499)$\\
\SetRow{lightgray!40} 1511 & 14 & $L_{3}(1511)$\\
1531 & 15 & $U_{3}(647)$, $U_{4}(647)$\\
\SetRow{lightgray!40} 1543 & 8 & $L_{3}(1543)$\\
1549 & 7 & $L_{2}(1549^{2})$, $S_{4}(1549)$\\
\SetRow{lightgray!40} 1553 & 9 & $L_{2}(1553^{2})$, $S_{4}(1553)$\\
1567 & 10 & $L_{3}(1031)$, $L_{2}(1031^{3})$, $G_2(1031)$, $L_{2}(1567^{2})$, $S_{4}(1567)$\\
\SetRow{lightgray!40} 1571 & 10 & $L_{3}(1571)$\\
1583 & 16 & $U_{3}(1583)$\\
\SetRow{lightgray!40} 1597 & 10 & $U_{3}(223)$, $L_{3}(1597)$, $L_{4}(1597)$, $L_{2}(1597^{2})$, $S_{4}(1597)$\\
1601 & 11 & $L_{3}(1601)$, $L_{2}(1601^{3})$, $G_2(1601)$, $U_{3}(1601)$\\
\SetRow{lightgray!40} 1607 & 5 & $L_{2}(1607^{2})$, $S_{4}(1607)$\\
1609 & 16 & $L_{3}(251^{2})$, $L_{2}(251^{3})$, $S_{6}(251)$, $O_{7}(251)$, $O^+_{8}(251)$, $G_2(251)$, $U_{3}(251)$, $U_{4}(251)$, $L_{2}(523^{2})$, $S_{4}(523)$, $U_{4}(523)$\\
\SetRow{lightgray!40} 1613 & 10 & $L_{2}(127^{2})$, $S_{4}(127)$, $U_{4}(127)$\\
1619 & 4 & $U_{3}(1619)$\\
\SetRow{lightgray!40} 1621 & 8 & $U_{3}(89^{2})$\\
1637 & 27 & $L_{2}(1321^{2})$, $S_{4}(1321)$, $L_{3}(1637)$, $L_{4}(1637)$, $L_{2}(1637^{2})$, $S_{4}(1637)$\\
\SetRow{lightgray!40} 1657 & 12 & $U_{3}(71)$, $L_{2}(239^{4})$, $S_{4}(239^{2})$, $L_{2}(1657^{2})$, $S_{4}(1657)$\\
1663 & 6 & $L_{3}(1663)$\\
\SetRow{lightgray!40} 1693 & 16 & $L_{3}(433)$, $L_{3}(1259)$, $L_{4}(1601)$, $L_{2}(1601^{2})$, $L_{3}(1601^{2})$, $S_{4}(1601)$, $S_{6}(1601)$, $O_{7}(1601)$, $O^+_{8}(1601)$, $U_{4}(1601)$, $L_{3}(1693)$\\
1697 & 6 & $L_{4}(1283)$, $L_{2}(1283^{2})$, $S_{4}(1283)$\\
\SetRow{lightgray!40} 1699 & 20 & $L_{3}(397)$, $L_{2}(397^{3})$, $G_2(397)$, $L_{3}(1301)$, $L_{4}(1301)$, $L_{2}(1699^{2})$, $S_{4}(1699)$, $U_{3}(1699)$, $U_{4}(1699)$\\
1709 & 15 & $L_{2}(1319^{2})$, $S_{4}(1319)$\\
\SetRow{lightgray!40} 1723 & 21 & $L_{3}(41)$, $L_{4}(41)$, $L_{3}(41^{2})$, $L_{2}(41^{3})$, $S_{6}(41)$, $O_{7}(41)$, $O^+_{8}(41)$, $G_2(41)$, $L_{2}(1723^{2})$, $S_{4}(1723)$\\
1741 & 14 & $L_{5}(5^{3})$, $L_{3}(5^{5})$, $L_{2}(59^{2})$, $L_{2}(59^{4})$, $S_{4}(59)$, $S_{4}(59^{2})$, $U_{4}(59)$\\
\SetRow{lightgray!40} 1747 & 9 & $L_{2}(1747^{2})$, $S_{4}(1747)$\\
1753 & 11 & $L_{2}(1571^{3})$, $G_2(1571)$, $U_{3}(1571)$, $U_{3}(1753)$\\
\SetRow{lightgray!40} 1759 & 21 & $U_{3}(509)$, $U_{4}(509)$\\
1783 & 19 & $L_{3}(193)$, $L_{4}(193)$, $L_{3}(1783)$, $L_{4}(1783)$, $L_{2}(1783^{2})$, $L_{3}(1783^{2})$, $L_{2}(1783^{3})$, $S_{4}(1783)$, $S_{6}(1783)$, $O_{7}(1783)$, $O^+_{8}(1783)$, $G_2(1783)$, $U_{3}(1783)$, $U_{4}(1783)$\\
\SetRow{lightgray!40} 1787 & 5 & $L_{2}(1787^{2})$, $S_{4}(1787)$\\
1789 & 25 & $L_{3}(1637^{2})$, $L_{2}(1637^{3})$, $S_{6}(1637)$, $O_{7}(1637)$, $O^+_{8}(1637)$, $G_2(1637)$, $U_{3}(1637)$, $U_{4}(1637)$, $L_{3}(1789)$, $L_{4}(1789)$, $L_{2}(1789^{2})$, $S_{4}(1789)$\\
\SetRow{lightgray!40} 1801 & 23 & $L_{5}(32)$, $L_{6}(32)$, $L_{3}(73)$, $L_{4}(73)$, $L_{3}(73^{2})$, $L_{2}(73^{3})$, $S_{6}(73)$, $O_{7}(73)$, $O^+_{8}(73)$, $G_2(73)$, $L_{2}(977^{2})$, $S_{4}(977)$\\
1823 & 21 & $L_{3}(1823)$, $L_{4}(1823)$, $L_{2}(1823^{2})$, $L_{3}(1823^{2})$, $L_{2}(1823^{3})$, $S_{4}(1823)$, $S_{6}(1823)$, $O_{7}(1823)$, $O^+_{8}(1823)$, $G_2(1823)$, $U_{3}(1823)$, $U_{4}(1823)$\\
\SetRow{lightgray!40} 1831 & 19 & $U_{3}(673)$, $L_{3}(1831)$\\
1861 & 17 & $L_{4}(61)$, $L_{2}(61^{2})$, $L_{3}(61^{2})$, $S_{4}(61)$, $S_{6}(61)$, $O_{7}(61)$, $O^+_{8}(61)$, $U_{4}(61)$, $L_{2}(1861^{2})$, $S_{4}(1861)$\\
\SetRow{lightgray!40} 1867 & 6 & $U_{3}(1033)$\\
1871 & 4 & $U_{3}(1871)$\\
\SetRow{lightgray!40} 1873 & 7 & $U_{3}(1759)$, $L_{3}(1873)$\\
1877 & 8 & $L_{4}(137)$, $L_{2}(137^{2})$, $S_{4}(137)$, $L_{2}(1877^{2})$, $S_{4}(1877)$\\
\SetRow{lightgray!40} 1889 & 17 & $L_{2}(331^{2})$, $S_{4}(331)$, $L_{2}(1889^{2})$, $S_{4}(1889)$\\
1913 & 25 & $L_{2}(1201^{2})$, $S_{4}(1201)$, $L_{3}(1913)$, $L_{4}(1913)$, $L_{2}(1913^{2})$, $S_{4}(1913)$\\
\SetRow{lightgray!40} 1931 & 7 & $L_{2}(1931^{2})$, $S_{4}(1931)$, $U_{3}(1931)$, $U_{4}(1931)$\\
1951 & 24 & $U_{3}(1951)$\\
\SetRow{lightgray!40} 1973 & 8 & $U_{3}(1973)$\\
1979 & 10 & $U_{3}(1979)$\\
\SetRow{lightgray!40} 1987 & 16 & $L_{3}(647)$, $L_{4}(647)$, $L_{3}(647^{2})$, $L_{2}(647^{3})$, $S_{6}(647)$, $O_{7}(647)$, $O^+_{8}(647)$, $G_2(647)$, $L_{3}(1987)$\\
1993 & 18 & $L_{2}(41^{6})$, $S_{4}(41^{3})$, $G_2(41^{2})$, $^3D_4(41)$, $U_{3}(41^{2})$, $L_{3}(313^{2})$, $L_{2}(313^{3})$, $S_{6}(313)$, $O_{7}(313)$, $O^+_{8}(313)$, $G_2(313)$, $U_{3}(313)$, $U_{4}(313)$\\
\SetRow{lightgray!40} 1999 & 13 & $L_{3}(809^{2})$, $L_{2}(809^{3})$, $S_{6}(809)$, $O_{7}(809)$, $O^+_{8}(809)$, $G_2(809)$, $U_{3}(809)$, $U_{4}(809)$\\
2017 & 16 & $L_{4}(229)$, $L_{2}(229^{2})$, $S_{4}(229)$, $U_{3}(1723)$, $U_{4}(1723)$\\
\SetRow{lightgray!40} 2027 & 4 & $L_{3}(2027)$\\
2039 & 16 & $L_{3}(2039)$\\
\SetRow{lightgray!40} 2053 & 15 & $L_{3}(197)$, $L_{2}(197^{3})$, $G_2(197)$, $L_{3}(2053)$\\
2081 & 6 & $L_{2}(1979^{2})$, $S_{4}(1979)$, $U_{4}(1979)$\\
\SetRow{lightgray!40} 2083 & 7 & $L_{3}(449)$, $U_{3}(2083)$\\
2089 & 12 & $U_{3}(827)$\\
\SetRow{lightgray!40} 2099 & 17 & $L_{2}(2099^{2})$, $S_{4}(2099)$, $U_{3}(2099)$, $U_{4}(2099)$\\
2111 & 5 & $L_{2}(2111^{2})$, $S_{4}(2111)$\\
\SetRow{lightgray!40} 2113 & 34 & $L_{2}(2^{22})$, $S_{4}(2^{11})$, $U_{11}(4)$, $U_{12}(4)$, $U_{13}(4)$, $U_{14}(4)$, $U_{15}(4)$, $L_{3}(439^{2})$, $L_{2}(439^{3})$, $S_{6}(439)$, $O_{7}(439)$, $O^+_{8}(439)$, $G_2(439)$, $U_{3}(439)$, $U_{4}(439)$, $L_{3}(2113)$, $Sz(2^{11})$\\
2131 & 10 & $L_{2}(1663^{3})$, $G_2(1663)$, $U_{3}(1663)$\\
\SetRow{lightgray!40} 2141 & 6 & $L_{2}(419^{2})$, $S_{4}(419)$, $U_{4}(419)$\\
2143 & 14 & $L_{3}(349)$, $L_{2}(2143^{2})$, $S_{4}(2143)$\\
\SetRow{lightgray!40} 2153 & 11 & $L_{2}(2153^{2})$, $S_{4}(2153)$\\
2161 & 29 & $L_{3}(593)$, $L_{4}(593)$, $L_{3}(593^{2})$, $L_{2}(593^{3})$, $S_{6}(593)$, $O_{7}(593)$, $O^+_{8}(593)$, $G_2(593)$, $L_{3}(1567)$, $L_{4}(1567)$\\
\SetRow{lightgray!40} 2203 & 6 & $U_{3}(2203)$\\
2207 & 9 & $L_{2}(2207^{2})$, $S_{4}(2207)$\\
\SetRow{lightgray!40} 2213 & 10 & $U_{3}(2213)$\\
2237 & 5 & $L_{2}(1021^{2})$, $S_{4}(1021)$\\
\SetRow{lightgray!40} 2243 & 10 & $L_{3}(2243)$\\
2251 & 24 & $L_{2}(19^{5})$, $U_{5}(19)$, $U_{6}(19)$, $U_{3}(709)$, $L_{2}(1543^{3})$, $G_2(1543)$, $U_{3}(1543)$\\
\SetRow{lightgray!40} 2267 & 7 & $L_{3}(2267)$, $L_{2}(2267^{3})$, $G_2(2267)$, $U_{3}(2267)$\\
2269 & 16 & $U_{7}(27)$, $U_{3}(3^{7})$, $L_{3}(83^{2})$, $L_{2}(83^{3})$, $S_{6}(83)$, $O_{7}(83)$, $O^+_{8}(83)$, $G_2(83)$, $U_{3}(83)$, $U_{4}(83)$, $^2G_2(3^{7})$\\
\SetRow{lightgray!40} 2273 & 10 & $U_{3}(2273)$\\
2281 & 15 & $L_{4}(1571)$, $L_{2}(1571^{2})$, $L_{3}(1571^{2})$, $S_{4}(1571)$, $S_{6}(1571)$, $O_{7}(1571)$, $O^+_{8}(1571)$, $U_{4}(1571)$\\
\SetRow{lightgray!40} 2287 & 9 & $U_{3}(1483)$, $U_{4}(1483)$\\
2293 & 12 & $U_{3}(113^{2})$, $L_{3}(1303)$, $L_{4}(1303)$, $L_{4}(1693)$, $L_{2}(1693^{2})$, $S_{4}(1693)$, $L_{3}(2293)$\\
\SetRow{lightgray!40} 2297 & 15 & $L_{2}(2297^{2})$, $S_{4}(2297)$\\
2309 & 7 & $L_{2}(1621^{2})$, $S_{4}(1621)$, $L_{2}(2309^{2})$, $S_{4}(2309)$\\
\SetRow{lightgray!40} 2311 & 28 & $U_{3}(883)$, $L_{2}(1429^{3})$, $G_2(1429)$, $U_{3}(1429)$, $L_{3}(2311)$\\
2333 & 9 & $L_{2}(2333^{2})$, $S_{4}(2333)$\\
\SetRow{lightgray!40} 2339 & 4 & $U_{3}(2339)$\\
2347 & 7 & $U_{3}(1063)$, $L_{3}(2347)$\\
\SetRow{lightgray!40} 2371 & 8 & $U_{3}(1907)$\\
2383 & 10 & $L_{3}(1103)$, $L_{3}(1279)$, $U_{3}(2383)$\\
\SetRow{lightgray!40} 2389 & 13 & $L_{3}(1699)$, $L_{4}(1699)$, $L_{3}(1699^{2})$, $L_{2}(1699^{3})$, $S_{6}(1699)$, $O_{7}(1699)$, $O^+_{8}(1699)$, $G_2(1699)$\\
2393 & 14 & $L_{4}(971)$, $L_{2}(971^{2})$, $S_{4}(971)$, $L_{3}(2393)$, $L_{2}(2393^{3})$, $G_2(2393)$, $U_{3}(2393)$\\
\SetRow{lightgray!40} 2411 & 10 & $U_{5}(13)$, $U_{6}(13)$, $U_{3}(2411)$\\
2417 & 9 & $L_{2}(2417^{2})$, $S_{4}(2417)$\\
\SetRow{lightgray!40} 2437 & 9 & $L_{4}(2039)$, $L_{2}(2039^{2})$, $S_{4}(2039)$, $L_{3}(2351)$\\
2441 & 9 & $L_{2}(2441^{2})$, $S_{4}(2441)$\\
\SetRow{lightgray!40} 2467 & 8 & $U_{3}(2251)$\\
2473 & 6 & $L_{3}(2473)$\\
\SetRow{lightgray!40} 2503 & 25 & $U_{3}(1277)$, $U_{4}(1277)$, $L_{3}(2503)$, $L_{4}(2503)$, $L_{2}(2503^{2})$, $S_{4}(2503)$\\
2521 & 15 & $L_{2}(71^{2})$, $S_{4}(71)$, $U_{4}(71)$, $U_{3}(97^{2})$\\
\SetRow{lightgray!40} 2531 & 11 & $L_{2}(2531^{2})$, $S_{4}(2531)$\\
2539 & 13 & $L_{3}(307^{2})$, $L_{2}(307^{3})$, $S_{6}(307)$, $O_{7}(307)$, $O^+_{8}(307)$, $G_2(307)$, $U_{3}(307)$, $U_{4}(307)$\\
\SetRow{lightgray!40} 2551 & 11 & $L_{3}(2551)$, $L_{4}(2551)$, $L_{2}(2551^{2})$, $S_{4}(2551)$\\
2557 & 25 & $L_{3}(1721)$, $L_{3}(2557)$\\
\SetRow{lightgray!40} 2579 & 17 & $L_{2}(2579^{2})$, $S_{4}(2579)$, $U_{3}(2579)$, $U_{4}(2579)$\\
2591 & 4 & $L_{3}(2591)$\\
\SetRow{lightgray!40} 2609 & 15 & $L_{2}(389^{2})$, $S_{4}(389)$, $L_{3}(2609)$, $L_{2}(2609^{3})$, $G_2(2609)$, $U_{3}(2609)$\\
2617 & 8 & $U_{3}(1553)$, $U_{4}(1553)$, $U_{3}(2617)$\\
\SetRow{lightgray!40} 2621 & 15 & $L_{2}(2621^{2})$, $S_{4}(2621)$\\
2633 & 19 & $L_{2}(1409^{2})$, $S_{4}(1409)$, $L_{2}(2633^{2})$, $S_{4}(2633)$\\
\SetRow{lightgray!40} 2657 & 7 & $L_{4}(163)$, $L_{2}(163^{2})$, $L_{5}(163)$, $S_{4}(163)$\\
2671 & 8 & $L_{3}(2671)$\\
\SetRow{lightgray!40} 2677 & 10 & $L_{3}(1033)$, $L_{2}(1033^{3})$, $G_2(1033)$\\
2683 & 6 & $L_{3}(2683)$\\
\SetRow{lightgray!40} 2687 & 4 & $U_{3}(2687)$\\
2689 & 7 & $L_{3}(2297)$, $L_{4}(2297)$\\
\SetRow{lightgray!40} 2693 & 10 & $L_{2}(859^{2})$, $S_{4}(859)$, $U_{4}(859)$\\
2707 & 15 & $L_{3}(1327)$, $L_{4}(1327)$, $L_{3}(1327^{2})$, $L_{2}(1327^{3})$, $S_{6}(1327)$, $O_{7}(1327)$, $O^+_{8}(1327)$, $G_2(1327)$, $L_{2}(2707^{2})$, $S_{4}(2707)$\\
\SetRow{lightgray!40} 2713 & 9 & $L_{2}(887^{2})$, $S_{4}(887)$\\
2719 & 12 & $L_{3}(1453)$\\
\SetRow{lightgray!40} 2729 & 6 & $L_{2}(1627^{2})$, $S_{4}(1627)$, $U_{3}(2729)$\\
2731 & 16 & $O^-_{26}(2)$, $U_{13}(2)$, $U_{14}(2)$, $U_{15}(2)$, $U_{16}(2)$\\
\SetRow{lightgray!40} 2749 & 7 & $L_{3}(2153)$, $L_{4}(2153)$\\
2753 & 17 & $L_{2}(2753^{2})$, $S_{4}(2753)$\\
\SetRow{lightgray!40} 2767 & 13 & $L_{2}(2767^{2})$, $S_{4}(2767)$\\
2777 & 15 & $L_{2}(2777^{2})$, $S_{4}(2777)$\\
\SetRow{lightgray!40} 2789 & 6 & $L_{2}(167^{2})$, $S_{4}(167)$, $U_{3}(2789)$\\
2791 & 8 & $L_{3}(2699)$\\
\SetRow{lightgray!40} 2797 & 6 & $U_{3}(1697)$\\
2801 & 27 & $L_{5}(7)$, $L_{6}(7)$, $L_{5}(49)$, $L_{2}(7^{5})$, $L_{6}(49)$, $S_{10}(7)$, $S_{12}(7)$, $O_{11}(7)$, $O_{13}(7)$, $O^+_{10}(7)$, $O^+_{12}(7)$, $O^-_{12}(7)$, $O^-_{14}(7)$, $E_6(7)$, $L_{4}(1543)$, $L_{2}(1543^{2})$, $L_{3}(1543^{2})$, $S_{4}(1543)$, $S_{6}(1543)$, $O_{7}(1543)$, $O^+_{8}(1543)$, $U_{4}(1543)$, $L_{2}(2801^{2})$, $S_{4}(2801)$\\
\SetRow{lightgray!40} 2803 & 18 & $L_{3}(2389)$\\
2819 & 16 & $L_{3}(2819)$\\
\SetRow{lightgray!40} 2833 & 14 & $L_{3}(1301^{2})$, $L_{2}(1301^{3})$, $S_{6}(1301)$, $O_{7}(1301)$, $O^+_{8}(1301)$, $G_2(1301)$, $U_{3}(1301)$, $U_{4}(1301)$, $U_{3}(2833)$\\
2837 & 8 & $L_{3}(2837)$\\
\SetRow{lightgray!40} 2843 & 13 & $L_{2}(2843^{2})$, $S_{4}(2843)$, $U_{3}(2843)$, $U_{4}(2843)$\\
2851 & 9 & $U_{5}(107)$, $U_{6}(107)$\\
\SetRow{lightgray!40} 2861 & 22 & $L_{5}(149)$, $L_{2}(2861^{2})$, $S_{4}(2861)$\\
2897 & 9 & $L_{2}(1777^{2})$, $S_{4}(1777)$\\
\SetRow{lightgray!40} 2903 & 9 & $L_{2}(2903^{2})$, $S_{4}(2903)$\\
2917 & 13 & $L_{2}(2917^{2})$, $S_{4}(2917)$\\
\SetRow{lightgray!40} 2939 & 19 & $L_{3}(2939)$, $L_{4}(2939)$, $L_{2}(2939^{2})$, $S_{4}(2939)$\\
2953 & 14 & $L_{3}(2153^{2})$, $L_{2}(2153^{3})$, $S_{6}(2153)$, $O_{7}(2153)$, $O^+_{8}(2153)$, $G_2(2153)$, $U_{3}(2153)$, $U_{4}(2153)$, $U_{3}(2953)$\\
\SetRow{lightgray!40} 2957 & 8 & $U_{3}(2957)$\\
2969 & 6 & $U_{7}(23)$, $L_{2}(2969^{2})$, $S_{4}(2969)$\\
\SetRow{lightgray!40} 2971 & 43 & $U_{3}(2917)$, $U_{4}(2917)$, $L_{3}(2971)$, $L_{4}(2971)$, $L_{2}(2971^{2})$, $L_{3}(2971^{2})$, $L_{2}(2971^{3})$, $S_{4}(2971)$, $S_{6}(2971)$, $O_{7}(2971)$, $O^+_{8}(2971)$, $G_2(2971)$, $U_{3}(2971)$, $U_{4}(2971)$\\
2999 & 4 & $L_{3}(2999)$\\
\SetRow{lightgray!40} 3001 & 15 & $L_{3}(3001)$, $L_{2}(3001^{3})$, $G_2(3001)$, $U_{3}(3001)$\\
3011 & 11 & $L_{2}(3011^{2})$, $S_{4}(3011)$\\
\SetRow{lightgray!40} 3019 & 18 & $L_{3}(239)$, $L_{4}(239)$, $L_{3}(239^{2})$, $L_{2}(239^{3})$, $L_{4}(239^{2})$, $S_{6}(239)$, $S_{8}(239)$, $O_{7}(239)$, $O_{9}(239)$, $O^+_{8}(239)$, $O^-_{8}(239)$, $G_2(239)$, $L_{3}(3019)$\\
3037 & 8 & $L_{4}(281)$, $L_{2}(281^{2})$, $S_{4}(281)$\\
\SetRow{lightgray!40} 3041 & 17 & $L_{4}(2267)$, $L_{2}(2267^{2})$, $L_{3}(2267^{2})$, $S_{4}(2267)$, $S_{6}(2267)$, $O_{7}(2267)$, $O^+_{8}(2267)$, $U_{4}(2267)$\\
3049 & 15 & $L_{2}(137^{4})$, $S_{4}(137^{2})$\\
\SetRow{lightgray!40} 3079 & 9 & $L_{3}(43^{3})$, $L_{2}(547^{3})$, $G_2(547)$, $U_{3}(547)$\\
3083 & 11 & $L_{2}(3083^{2})$, $S_{4}(3083)$, $U_{3}(3083)$, $U_{4}(3083)$\\
\SetRow{lightgray!40} 3089 & 22 & $L_{3}(3089)$\\
3109 & 13 & $L_{2}(727^{2})$, $S_{4}(727)$\\
\SetRow{lightgray!40} 3121 & 23 & $L_{4}(79)$, $L_{2}(79^{2})$, $L_{2}(79^{4})$, $S_{4}(79)$, $S_{4}(79^{2})$, $L_{3}(1999)$\\
3137 & 31 & $L_{2}(3137^{2})$, $S_{4}(3137)$, $U_{3}(3137)$, $U_{4}(3137)$\\
\SetRow{lightgray!40} 3169 & 25 & $L_{3}(97)$, $L_{4}(97)$, $L_{3}(97^{2})$, $L_{2}(97^{3})$, $L_{2}(97^{6})$, $S_{6}(97)$, $S_{4}(97^{3})$, $O_{7}(97)$, $O^+_{8}(97)$, $G_2(97)$, $G_2(97^{2})$, $^3D_4(97)$\\
3181 & 8 & $U_{3}(2741)$\\
\SetRow{lightgray!40} 3187 & 8 & $L_{3}(1871)$, $L_{2}(1871^{3})$, $G_2(1871)$\\
3203 & 11 & $L_{3}(3203)$, $L_{4}(3203)$, $L_{2}(3203^{2})$, $S_{4}(3203)$\\
\SetRow{lightgray!40} 3221 & 11 & $L_{5}(11)$, $L_{6}(11)$\\
3229 & 26 & $L_{2}(839^{2})$, $S_{4}(839)$, $L_{3}(3229)$\\
\SetRow{lightgray!40} 3251 & 5 & $L_{2}(3251^{2})$, $S_{4}(3251)$\\
3253 & 8 & $L_{3}(1439)$, $L_{2}(3253^{2})$, $S_{4}(3253)$\\
\SetRow{lightgray!40} 3257 & 4 & $L_{3}(3257)$\\
3259 & 15 & $U_{3}(853)$, $U_{4}(853)$\\
\SetRow{lightgray!40} 3301 & 10 & $L_{2}(2089^{2})$, $S_{4}(2089)$, $U_{3}(3011^{2})$\\
3307 & 9 & $L_{2}(3307^{2})$, $S_{4}(3307)$\\
\SetRow{lightgray!40} 3313 & 9 & $L_{3}(1123)$, $L_{4}(1123)$\\
3319 & 9 & $L_{3}(3319)$, $L_{4}(3319)$, $L_{2}(3319^{2})$, $S_{4}(3319)$\\
\SetRow{lightgray!40} 3323 & 11 & $L_{2}(3323^{2})$, $S_{4}(3323)$, $U_{3}(3323)$, $U_{4}(3323)$\\
3331 & 14 & $L_{3}(1867)$\\
\SetRow{lightgray!40} 3343 & 7 & $L_{2}(3343^{2})$, $S_{4}(3343)$\\
3347 & 14 & $L_{3}(3347)$\\
\SetRow{lightgray!40} 3359 & 7 & $L_{3}(3359)$, $L_{4}(3359)$, $L_{2}(3359^{2})$, $S_{4}(3359)$\\
3361 & 14 & $U_{3}(421^{2})$, $L_{2}(3361^{2})$, $S_{4}(3361)$\\
\SetRow{lightgray!40} 3373 & 22 & $U_{3}(2719)$, $L_{3}(3373)$, $L_{4}(3373)$, $L_{2}(3373^{2})$, $S_{4}(3373)$\\
3389 & 7 & $L_{3}(3389)$, $L_{4}(3389)$, $L_{2}(3389^{2})$, $S_{4}(3389)$\\
\SetRow{lightgray!40} 3391 & 18 & $U_{3}(3391)$\\
3413 & 24 & $L_{4}(1471)$, $L_{2}(1471^{2})$, $S_{4}(1471)$\\
\SetRow{lightgray!40} 3433 & 25 & $L_{3}(269^{2})$, $L_{2}(269^{3})$, $S_{6}(269)$, $O_{7}(269)$, $O^+_{8}(269)$, $G_2(269)$, $U_{3}(269)$, $U_{4}(269)$\\
3449 & 10 & $U_{3}(3449)$\\
\SetRow{lightgray!40} 3457 & 7 & $L_{2}(2749^{2})$, $S_{4}(2749)$\\
3461 & 7 & $L_{4}(1453)$, $L_{2}(1453^{2})$, $S_{4}(1453)$, $L_{3}(3461)$\\
\SetRow{lightgray!40} 3463 & 10 & $L_{3}(367)$, $L_{2}(367^{3})$, $G_2(367)$, $L_{2}(3463^{2})$, $S_{4}(3463)$\\
3467 & 4 & $U_{3}(3467)$\\
\SetRow{lightgray!40} 3499 & 17 & $U_{3}(157)$, $U_{4}(157)$, $U_{3}(3343)$, $U_{4}(3343)$\\
3511 & 11 & $U_{3}(757)$, $U_{4}(757)$, $L_{2}(3511^{2})$, $S_{4}(3511)$\\
\SetRow{lightgray!40} 3517 & 12 & $U_{3}(3259)$\\
3527 & 7 & $L_{3}(3527)$, $L_{2}(3527^{3})$, $G_2(3527)$, $U_{3}(3527)$\\
\SetRow{lightgray!40} 3529 & 9 & $L_{2}(449^{3})$, $G_2(449)$, $U_{3}(449)$, $U_{3}(3529)$\\
3533 & 8 & $U_{3}(3533)$\\
\SetRow{lightgray!40} 3541 & 29 & $L_{3}(59)$, $L_{4}(59)$, $L_{5}(59)$, $L_{6}(59)$, $L_{3}(59^{2})$, $L_{2}(59^{3})$, $L_{4}(59^{2})$, $S_{6}(59)$, $S_{8}(59)$, $O_{7}(59)$, $O_{9}(59)$, $O^+_{8}(59)$, $O^+_{10}(59)$, $O^-_{8}(59)$, $G_2(59)$, $U_{3}(59^{3})$, $L_{2}(2689^{2})$, $S_{4}(2689)$, $L_{3}(3541)$, $L_{4}(3541)$, $L_{2}(3541^{2})$, $S_{4}(3541)$\\
3547 & 14 & $U_{3}(1163)$, $L_{2}(3547^{2})$, $S_{4}(3547)$\\
\SetRow{lightgray!40} 3557 & 5 & $L_{2}(3557^{2})$, $S_{4}(3557)$\\
3559 & 14 & $U_{3}(3559)$\\
\SetRow{lightgray!40} 3571 & 23 & $L_{3}(103)$, $L_{4}(103)$, $L_{3}(103^{2})$, $L_{2}(103^{3})$, $S_{6}(103)$, $O_{7}(103)$, $O^+_{8}(103)$, $G_2(103)$, $L_{3}(3467)$, $L_{2}(3467^{3})$, $G_2(3467)$, $L_{3}(3571)$\\
3581 & 6 & $L_{2}(3217^{2})$, $S_{4}(3217)$, $L_{3}(3581)$\\
\SetRow{lightgray!40} 3583 & 14 & $U_{3}(1039)$, $L_{2}(3583^{2})$, $S_{4}(3583)$\\
3593 & 18 & $L_{4}(1153)$, $L_{2}(1153^{2})$, $S_{4}(1153)$\\
\SetRow{lightgray!40} 3607 & 13 & $L_{3}(1399)$, $L_{2}(1399^{3})$, $G_2(1399)$, $L_{3}(2207)$, $L_{4}(2207)$, $L_{3}(3607)$\\
3613 & 6 & $U_{3}(3613)$\\
\SetRow{lightgray!40} 3617 & 12 & $L_{2}(2383^{2})$, $S_{4}(2383)$, $U_{4}(2383)$, $L_{2}(3617^{2})$, $S_{4}(3617)$\\
3623 & 13 & $L_{3}(3623)$, $L_{2}(3623^{3})$, $G_2(3623)$, $U_{3}(3623)$\\
\SetRow{lightgray!40} 3631 & 12 & $U_{5}(523)$, $L_{3}(3631)$, $L_{2}(3631^{3})$, $G_2(3631)$, $U_{3}(3631)$\\
3637 & 8 & $U_{3}(3637)$\\
\SetRow{lightgray!40} 3643 & 20 & $U_{3}(3221)$, $L_{2}(3643^{2})$, $S_{4}(3643)$\\
3659 & 15 & $L_{2}(3659^{2})$, $S_{4}(3659)$\\
\SetRow{lightgray!40} 3671 & 5 & $L_{2}(3671^{2})$, $S_{4}(3671)$\\
3673 & 8 & $L_{3}(1151)$, $L_{3}(2521)$, $L_{3}(3673)$\\
\SetRow{lightgray!40} 3691 & 9 & $U_{3}(3217)$, $U_{4}(3217)$\\
3701 & 13 & $L_{4}(1279)$, $L_{2}(1279^{2})$, $S_{4}(1279)$, $L_{3}(3701)$\\
\SetRow{lightgray!40} 3709 & 22 & $L_{3}(499^{2})$, $L_{2}(499^{3})$, $S_{6}(499)$, $O_{7}(499)$, $O^+_{8}(499)$, $G_2(499)$, $U_{3}(499)$, $U_{4}(499)$, $L_{2}(1609^{2})$, $S_{4}(1609)$, $U_{3}(3709)$\\
3727 & 10 & $U_{3}(2539)$, $L_{2}(3727^{2})$, $S_{4}(3727)$\\
\SetRow{lightgray!40} 3739 & 24 & $L_{3}(3739)$\\
3769 & 31 & $L_{3}(463)$, $L_{4}(463)$, $L_{3}(463^{2})$, $L_{2}(463^{3})$, $S_{6}(463)$, $O_{7}(463)$, $O^+_{8}(463)$, $G_2(463)$, $L_{3}(3769)$, $L_{4}(3769)$, $L_{2}(3769^{2})$, $L_{3}(3769^{2})$, $L_{2}(3769^{3})$, $S_{4}(3769)$, $S_{6}(3769)$, $O_{7}(3769)$, $O^+_{8}(3769)$, $G_2(3769)$, $U_{3}(3769)$, $U_{4}(3769)$\\
\SetRow{lightgray!40} 3779 & 16 & $L_{3}(3779)$\\
3793 & 12 & $L_{2}(1069^{3})$, $G_2(1069)$, $U_{3}(1069)$, $L_{3}(3793)$, $L_{4}(3793)$, $L_{2}(3793^{2})$, $S_{4}(3793)$\\
\SetRow{lightgray!40} 3833 & 37 & $L_{4}(19^{2})$, $L_{2}(19^{4})$, $L_{5}(19^{2})$, $L_{6}(19^{2})$, $L_{3}(19^{4})$, $S_{8}(19)$, $S_{4}(19^{2})$, $S_{10}(19)$, $S_{12}(19)$, $S_{6}(19^{2})$, $O_{9}(19)$, $O_{11}(19)$, $O_{13}(19)$, $O_{7}(19^{2})$, $O^+_{10}(19)$, $O^+_{12}(19)$, $O^+_{8}(19^{2})$, $O^-_{8}(19)$, $O^-_{10}(19)$, $O^-_{12}(19)$, $F_4(19)$, $U_{4}(19^{2})$\\
3851 & 13 & $L_{11}(3)$, $L_{12}(3)$, $L_{11}(9)$, $L_{2}(3^{11})$, $S_{22}(3)$, $O_{23}(3)$, $O^+_{22}(3)$, $O^+_{24}(3)$, $U_{5}(53)$, $U_{6}(53)$\\
\SetRow{lightgray!40} 3853 & 13 & $L_{3}(2713)$, $U_{3}(3853)$\\
3863 & 17 & $L_{2}(3863^{2})$, $S_{4}(3863)$\\
\SetRow{lightgray!40} 3881 & 19 & $L_{4}(197)$, $L_{2}(197^{2})$, $L_{5}(197)$, $L_{6}(197)$, $L_{3}(197^{2})$, $S_{4}(197)$, $S_{6}(197)$, $O_{7}(197)$, $O^+_{8}(197)$, $U_{4}(197)$\\
3889 & 22 & $L_{2}(1999^{3})$, $G_2(1999)$, $U_{3}(1999)$\\
\SetRow{lightgray!40} 3907 & 7 & $L_{2}(3907^{2})$, $S_{4}(3907)$\\
3917 & 7 & $L_{3}(3917)$, $L_{2}(3917^{3})$, $G_2(3917)$, $U_{3}(3917)$\\
\SetRow{lightgray!40} 3919 & 7 & $L_{3}(2749)$, $L_{4}(2749)$\\
3931 & 15 & $L_{3}(617)$, $L_{3}(3313)$\\
\SetRow{lightgray!40} 3947 & 25 & $L_{2}(3947^{2})$, $S_{4}(3947)$, $U_{3}(3947)$, $U_{4}(3947)$\\
3967 & 25 & $U_{3}(3079)$, $U_{3}(3967)$\\
\SetRow{lightgray!40} 4003 & 9 & $L_{2}(823^{3})$, $G_2(823)$, $U_{3}(823)$, $U_{3}(3181)$\\
4007 & 11 & $L_{2}(4007^{2})$, $S_{4}(4007)$, $U_{3}(4007)$, $U_{4}(4007)$\\
\SetRow{lightgray!40} 4013 & 8 & $U_{3}(4013)$\\
4021 & 16 & $L_{2}(7^{10})$, $S_{4}(7^{5})$, $U_{5}(49)$, $U_{6}(49)$, $L_{2}(47^{6})$, $S_{4}(47^{3})$, $G_2(47^{2})$, $^3D_4(47)$, $U_{3}(47^{2})$\\
\SetRow{lightgray!40} 4027 & 43 & $L_{3}(2207^{2})$, $L_{2}(2207^{3})$, $S_{6}(2207)$, $O_{7}(2207)$, $O^+_{8}(2207)$, $G_2(2207)$, $U_{3}(2207)$, $U_{4}(2207)$, $L_{3}(4027)$, $L_{4}(4027)$, $L_{2}(4027^{2})$, $L_{3}(4027^{2})$, $L_{2}(4027^{3})$, $S_{4}(4027)$, $S_{6}(4027)$, $O_{7}(4027)$, $O^+_{8}(4027)$, $G_2(4027)$, $U_{3}(4027)$, $U_{4}(4027)$\\
4049 & 4 & $U_{3}(4049)$\\
\SetRow{lightgray!40} 4051 & 13 & $L_{2}(2^{25})$, $U_{5}(32)$, $U_{6}(32)$, $L_{3}(797)$, $L_{3}(3253)$, $L_{4}(3253)$\\
4057 & 19 & $U_{3}(1409)$, $U_{4}(1409)$\\
\SetRow{lightgray!40} 4073 & 11 & $L_{3}(4073)$, $L_{2}(4073^{3})$, $G_2(4073)$, $U_{3}(4073)$\\
4079 & 14 & $U_{3}(4079)$\\
\SetRow{lightgray!40} 4091 & 4 & $U_{3}(4091)$\\
4093 & 8 & $U_{3}(3191)$\\
\SetRow{lightgray!40} 4099 & 15 & $L_{3}(2017)$, $L_{3}(2081)$\\
4111 & 20 & $U_{5}(41)$, $U_{6}(41)$, $L_{3}(4111)$\\
\SetRow{lightgray!40} 4129 & 14 & $L_{3}(1979)$, $L_{4}(1979)$, $L_{3}(1979^{2})$, $L_{2}(1979^{3})$, $S_{6}(1979)$, $O_{7}(1979)$, $O^+_{8}(1979)$, $G_2(1979)$, $U_{3}(4129)$\\
4133 & 10 & $L_{2}(733^{2})$, $S_{4}(733)$, $L_{3}(4133)$\\
\SetRow{lightgray!40} 4177 & 31 & $L_{2}(457^{2})$, $S_{4}(457)$, $L_{2}(1103^{3})$, $G_2(1103)$, $U_{3}(1103)$, $U_{3}(4177)$\\
4201 & 12 & $U_{3}(1013^{2})$\\
\SetRow{lightgray!40} 4211 & 8 & $U_{3}(4211)$\\
4217 & 7 & $L_{2}(4217^{2})$, $S_{4}(4217)$, $U_{3}(4217)$, $U_{4}(4217)$\\
\SetRow{lightgray!40} 4219 & 23 & $L_{3}(113^{2})$, $L_{2}(113^{3})$, $L_{2}(113^{6})$, $S_{6}(113)$, $S_{4}(113^{3})$, $O_{7}(113)$, $O^+_{8}(113)$, $G_2(113)$, $G_2(113^{2})$, $^3D_4(113)$, $U_{3}(113)$, $U_{4}(113)$\\
4229 & 5 & $L_{2}(4229^{2})$, $S_{4}(4229)$\\
\SetRow{lightgray!40} 4243 & 12 & $U_{3}(4243)$\\
4253 & 11 & $L_{3}(4253)$, $L_{2}(4253^{3})$, $G_2(4253)$, $U_{3}(4253)$\\
\SetRow{lightgray!40} 4261 & 15 & $L_{2}(4261^{2})$, $S_{4}(4261)$, $U_{3}(4261)$, $U_{4}(4261)$\\
4271 & 7 & $L_{5}(37)$, $L_{6}(37)$, $U_{5}(599)$, $L_{3}(4271)$\\
\SetRow{lightgray!40} 4273 & 12 & $U_{3}(2663)$\\
4283 & 8 & $L_{3}(4283)$\\
\SetRow{lightgray!40} 4289 & 12 & $L_{2}(3761^{2})$, $S_{4}(3761)$, $L_{3}(4289)$\\
4297 & 33 & $U_{3}(2887)$, $U_{3}(4297)$\\
\SetRow{lightgray!40} 4327 & 13 & $L_{2}(4327^{2})$, $S_{4}(4327)$\\
4337 & 5 & $L_{2}(4337^{2})$, $S_{4}(4337)$\\
\SetRow{lightgray!40} 4357 & 11 & $U_{3}(1319)$, $U_{4}(1319)$, $L_{2}(4357^{2})$, $S_{4}(4357)$\\
4363 & 15 & $L_{2}(4363^{2})$, $S_{4}(4363)$, $U_{3}(4363)$, $U_{4}(4363)$\\
\SetRow{lightgray!40} 4373 & 21 & $L_{2}(4373^{2})$, $S_{4}(4373)$\\
4391 & 9 & $L_{2}(4391^{2})$, $S_{4}(4391)$\\
\SetRow{lightgray!40} 4421 & 5 & $L_{2}(3469^{2})$, $S_{4}(3469)$\\
4423 & 31 & $L_{3}(67^{2})$, $L_{2}(67^{3})$, $S_{6}(67)$, $O_{7}(67)$, $O^+_{8}(67)$, $G_2(67)$, $U_{3}(67)$, $U_{4}(67)$, $U_{3}(4357)$, $U_{4}(4357)$, $L_{2}(4423^{2})$, $S_{4}(4423)$\\
\SetRow{lightgray!40} 4441 & 8 & $L_{3}(3539)$\\
4457 & 8 & $L_{3}(4457)$\\
\SetRow{lightgray!40} 4481 & 4 & $L_{3}(4481)$\\
4483 & 13 & $L_{2}(4483^{2})$, $S_{4}(4483)$\\
\SetRow{lightgray!40} 4493 & 20 & $L_{2}(2213^{2})$, $S_{4}(2213)$, $U_{4}(2213)$, $L_{2}(4493^{2})$, $S_{4}(4493)$\\
4513 & 7 & $L_{2}(4513^{2})$, $S_{4}(4513)$\\
\SetRow{lightgray!40} 4519 & 7 & $U_{3}(3463)$, $U_{4}(3463)$\\
4523 & 26 & $L_{3}(4523)$\\
\SetRow{lightgray!40} 4549 & 14 & $U_{3}(4549)$\\
4561 & 18 & $L_{5}(27)$, $L_{3}(3^{5})$, $L_{6}(27)$, $L_{4}(3^{5})$, $L_{3}(3^{10})$, $L_{2}(3^{15})$, $S_{6}(3^{5})$, $O_{7}(3^{5})$, $O^+_{8}(3^{5})$, $G_2(3^{5})$, $L_{3}(4561)$\\
\SetRow{lightgray!40} 4567 & 19 & $L_{2}(4567^{2})$, $S_{4}(4567)$\\
4591 & 9 & $U_{3}(311)$, $U_{4}(311)$\\
\SetRow{lightgray!40} 4597 & 10 & $L_{2}(2129^{2})$, $S_{4}(2129)$, $L_{3}(4219)$\\
4603 & 29 & $L_{3}(179)$, $L_{4}(179)$, $L_{3}(179^{2})$, $L_{2}(179^{3})$, $S_{6}(179)$, $O_{7}(179)$, $O^+_{8}(179)$, $G_2(179)$, $L_{3}(4423)$, $L_{4}(4423)$\\
\SetRow{lightgray!40} 4621 & 22 & $L_{3}(2857)$, $L_{3}(4621)$, $L_{2}(4621^{3})$, $G_2(4621)$, $U_{3}(4621)$\\
4637 & 6 & $L_{2}(2593^{2})$, $S_{4}(2593)$, $L_{3}(4637)$\\
\SetRow{lightgray!40} 4639 & 6 & $U_{3}(1361)$\\
4643 & 11 & $L_{3}(4643)$, $L_{4}(4643)$, $L_{2}(4643^{2})$, $S_{4}(4643)$\\
\SetRow{lightgray!40} 4649 & 7 & $L_{2}(2803^{2})$, $S_{4}(2803)$, $L_{2}(4649^{2})$, $S_{4}(4649)$\\
4651 & 16 & $L_{3}(787^{2})$, $L_{2}(787^{3})$, $S_{6}(787)$, $O_{7}(787)$, $O^+_{8}(787)$, $G_2(787)$, $U_{3}(787)$, $U_{4}(787)$, $L_{3}(4651)$\\
\SetRow{lightgray!40} 4657 & 9 & $L_{3}(967)$, $L_{3}(4657)$\\
4673 & 9 & $L_{2}(1993^{2})$, $S_{4}(1993)$\\
\SetRow{lightgray!40} 4679 & 17 & $L_{2}(4679^{2})$, $S_{4}(4679)$, $U_{3}(4679)$, $U_{4}(4679)$\\
4691 & 16 & $L_{7}(59)$, $L_{8}(59)$, $L_{3}(4691)$\\
\SetRow{lightgray!40} 4721 & 7 & $L_{2}(1697^{2})$, $S_{4}(1697)$, $U_{4}(1697)$, $U_{3}(1697^{2})$\\
4723 & 8 & $U_{3}(4723)$\\
\SetRow{lightgray!40} 4729 & 7 & $U_{3}(2693)$, $L_{3}(4729)$\\
4733 & 31 & $L_{7}(7)$, $L_{8}(7)$, $L_{9}(7)$, $L_{10}(7)$, $L_{7}(49)$, $L_{2}(7^{7})$, $S_{14}(7)$, $O_{15}(7)$, $O^+_{14}(7)$, $O^+_{16}(7)$, $L_{2}(4733^{2})$, $S_{4}(4733)$\\
\SetRow{lightgray!40} 4751 & 11 & $L_{2}(4751^{2})$, $S_{4}(4751)$\\
4783 & 7 & $L_{3}(3037)$, $L_{3}(4783)$\\
\SetRow{lightgray!40} 4787 & 5 & $L_{2}(4787^{2})$, $S_{4}(4787)$\\
4789 & 12 & $L_{2}(1481^{2})$, $S_{4}(1481)$, $L_{3}(3109)$, $L_{2}(4789^{2})$, $S_{4}(4789)$, $U_{3}(4789)$, $U_{4}(4789)$\\
\SetRow{lightgray!40} 4793 & 10 & $L_{4}(3313)$, $L_{2}(3313^{2})$, $S_{4}(3313)$\\
4799 & 4 & $L_{3}(4799)$\\
\SetRow{lightgray!40} 4801 & 14 & $U_{3}(2341)$\\
4813 & 7 & $U_{3}(1889)$, $U_{4}(1889)$\\
\SetRow{lightgray!40} 4817 & 20 & $L_{4}(1291)$, $L_{2}(1291^{2})$, $S_{4}(1291)$, $L_{2}(4817^{2})$, $S_{4}(4817)$\\
4831 & 33 & $L_{2}(4831^{2})$, $S_{4}(4831)$\\
\SetRow{lightgray!40} 4861 & 12 & $L_{3}(4861)$\\
4877 & 21 & $L_{4}(719)$, $L_{2}(719^{2})$, $L_{3}(719^{2})$, $S_{4}(719)$, $S_{6}(719)$, $O_{7}(719)$, $O^+_{8}(719)$, $U_{4}(719)$\\
\SetRow{lightgray!40} 4889 & 17 & $L_{2}(4159^{2})$, $S_{4}(4159)$\\
4903 & 13 & $U_{3}(2417)$, $U_{4}(2417)$, $L_{3}(4903)$, $L_{2}(4903^{3})$, $G_2(4903)$, $U_{3}(4903)$\\
\SetRow{lightgray!40} 4909 & 14 & $L_{2}(1613^{2})$, $S_{4}(1613)$, $L_{3}(4909)$\\
4933 & 15 & $L_{3}(2131)$, $L_{3}(2801)$, $L_{4}(2801)$, $L_{4}(3739)$, $L_{2}(3739^{2})$, $S_{4}(3739)$, $L_{3}(4933)$, $L_{4}(4933)$, $L_{2}(4933^{2})$, $S_{4}(4933)$\\
\SetRow{lightgray!40} 4937 & 9 & $L_{2}(4937^{2})$, $S_{4}(4937)$\\
4951 & 13 & $L_{3}(2689)$, $L_{4}(2689)$, $L_{3}(4951)$, $L_{2}(4951^{3})$, $G_2(4951)$, $U_{3}(4951)$\\
\SetRow{lightgray!40} 4957 & 16 & $L_{4}(359)$, $L_{2}(359^{2})$, $S_{4}(359)$, $L_{2}(4957^{2})$, $S_{4}(4957)$\\
4967 & 7 & $L_{2}(4967^{2})$, $S_{4}(4967)$, $U_{3}(4967)$, $U_{4}(4967)$\\
\SetRow{lightgray!40} 4969 & 12 & $L_{2}(4783^{3})$, $G_2(4783)$, $U_{3}(4783)$, $L_{3}(4969)$, $L_{4}(4969)$, $L_{2}(4969^{2})$, $S_{4}(4969)$\\
4973 & 19 & $L_{2}(223^{2})$, $S_{4}(223)$, $U_{4}(223)$, $U_{3}(4973)$\\
\SetRow{lightgray!40} 4987 & 12 & $U_{3}(3851)$, $L_{2}(4987^{2})$, $S_{4}(4987)$, $U_{3}(4987)$, $U_{4}(4987)$\\
5009 & 4 & $U_{3}(5009)$\\
\SetRow{lightgray!40} 5023 & 22 & $L_{3}(953)$, $L_{2}(953^{3})$, $G_2(953)$, $L_{2}(5023^{2})$, $S_{4}(5023)$\\
5051 & 10 & $L_{3}(5051)$\\
\SetRow{lightgray!40} 5059 & 28 & $L_{3}(1913^{2})$, $L_{2}(1913^{3})$, $S_{6}(1913)$, $O_{7}(1913)$, $O^+_{8}(1913)$, $G_2(1913)$, $U_{3}(1913)$, $U_{4}(1913)$, $U_{3}(5059)$\\
5077 & 12 & $L_{4}(4219)$, $L_{2}(4219^{2})$, $S_{4}(4219)$, $L_{3}(5077)$, $L_{2}(5077^{3})$, $G_2(5077)$, $U_{3}(5077)$\\
\SetRow{lightgray!40} 5087 & 15 & $L_{2}(5087^{2})$, $S_{4}(5087)$\\
5099 & 7 & $L_{3}(5099)$, $L_{4}(5099)$, $L_{2}(5099^{2})$, $S_{4}(5099)$\\
\SetRow{lightgray!40} 5101 & 11 & $L_{2}(101^{2})$, $S_{4}(101)$, $U_{4}(101)$, $L_{3}(5101)$\\
5107 & 15 & $L_{3}(311)$, $L_{4}(311)$, $L_{3}(311^{2})$, $L_{2}(311^{3})$, $S_{6}(311)$, $O_{7}(311)$, $O^+_{8}(311)$, $G_2(311)$\\
\SetRow{lightgray!40} 5113 & 21 & $L_{3}(71)$, $L_{4}(71)$, $L_{5}(71)$, $L_{6}(71)$, $L_{3}(71^{2})$, $L_{2}(71^{3})$, $S_{6}(71)$, $O_{7}(71)$, $O^+_{8}(71)$, $G_2(71)$, $L_{2}(5113^{2})$, $S_{4}(5113)$, $U_{3}(5113)$, $U_{4}(5113)$\\
5119 & 33 & $L_{2}(5119^{2})$, $S_{4}(5119)$, $U_{3}(5119)$, $U_{4}(5119)$\\
\SetRow{lightgray!40} 5153 & 27 & $L_{4}(227)$, $L_{2}(227^{2})$, $L_{3}(227^{2})$, $S_{4}(227)$, $S_{6}(227)$, $O_{7}(227)$, $O^+_{8}(227)$, $U_{4}(227)$, $L_{2}(5153^{2})$, $S_{4}(5153)$, $U_{3}(5153)$, $U_{4}(5153)$\\
5167 & 18 & $L_{6}(5^{3})$, $L_{3}(5^{6})$, $L_{2}(5^{9})$, $S_{6}(5^{3})$, $O_{7}(5^{3})$, $O^+_{8}(5^{3})$, $G_2(5^{3})$, $^2E_6(5)$, $U_{9}(5)$, $U_{3}(5^{3})$, $U_{10}(5)$, $U_{4}(5^{3})$, $U_{3}(5167)$\\
\SetRow{lightgray!40} 5171 & 11 & $L_{2}(5171^{2})$, $S_{4}(5171)$\\
5189 & 10 & $U_{3}(5189)$\\
\SetRow{lightgray!40} 5197 & 23 & $U_{3}(1879)$, $L_{3}(3319^{2})$, $L_{2}(3319^{3})$, $S_{6}(3319)$, $O_{7}(3319)$, $O^+_{8}(3319)$, $G_2(3319)$, $U_{3}(3319)$, $U_{4}(3319)$, $U_{3}(5197)$\\
5209 & 20 & $U_{3}(1193)$\\
\SetRow{lightgray!40} 5227 & 7 & $L_{2}(5227^{2})$, $S_{4}(5227)$\\
5231 & 5 & $L_{5}(307)$, $L_{6}(307)$\\
\SetRow{lightgray!40} 5233 & 7 & $L_{3}(331)$, $L_{4}(331)$\\
5237 & 26 & $U_{3}(5237)$\\
\SetRow{lightgray!40} 5261 & 17 & $L_{2}(827^{2})$, $S_{4}(827)$, $U_{4}(827)$, $L_{3}(5261)$\\
5281 & 21 & $U_{11}(5)$, $U_{12}(5)$, $L_{3}(3877)$, $U_{3}(5281)$\\
\SetRow{lightgray!40} 5323 & 21 & $L_{3}(1283^{2})$, $L_{2}(1283^{3})$, $S_{6}(1283)$, $O_{7}(1283)$, $O^+_{8}(1283)$, $G_2(1283)$, $U_{3}(1283)$, $U_{4}(1283)$, $L_{2}(5323^{2})$, $S_{4}(5323)$\\
5347 & 9 & $L_{3}(479)$, $L_{4}(479)$, $L_{2}(5347^{2})$, $S_{4}(5347)$\\
\SetRow{lightgray!40} 5351 & 32 & $U_{3}(5351)$\\
5381 & 8 & $L_{3}(5381)$\\
\SetRow{lightgray!40} 5413 & 7 & $L_{2}(5413^{2})$, $S_{4}(5413)$\\
5419 & 38 & $L_{7}(64)$, $L_{3}(2^{14})$, $L_{2}(2^{21})$, $L_{8}(64)$, $S_{14}(8)$, $S_{6}(2^{7})$, $S_{16}(8)$, $O^+_{16}(8)$, $O^+_{8}(2^{7})$, $O^-_{14}(8)$, $O^-_{16}(8)$, $G_2(2^{7})$, $U_{7}(8)$, $U_{3}(2^{7})$, $U_{8}(8)$, $U_{4}(2^{7})$, $L_{3}(127)$, $L_{4}(127)$, $L_{3}(127^{2})$, $L_{2}(127^{3})$, $S_{6}(127)$, $O_{7}(127)$, $O^+_{8}(127)$, $G_2(127)$, $L_{3}(5419)$\\
\SetRow{lightgray!40} 5431 & 8 & $U_{3}(5431)$\\
5437 & 6 & $L_{3}(5437)$\\
\SetRow{lightgray!40} 5441 & 7 & $L_{3}(5441)$, $L_{2}(5441^{3})$, $G_2(5441)$, $U_{3}(5441)$\\
5443 & 9 & $L_{2}(5443^{2})$, $S_{4}(5443)$\\
\SetRow{lightgray!40} 5477 & 7 & $L_{2}(5477^{2})$, $S_{4}(5477)$, $U_{3}(5477)$, $U_{4}(5477)$\\
5479 & 11 & $U_{3}(2777)$, $U_{4}(2777)$, $L_{3}(5479)$, $L_{4}(5479)$, $L_{2}(5479^{2})$, $S_{4}(5479)$\\
\SetRow{lightgray!40} 5483 & 23 & $L_{3}(5483)$, $L_{2}(5483^{3})$, $G_2(5483)$, $U_{3}(5483)$\\
5503 & 12 & $L_{3}(929)$, $L_{2}(929^{3})$, $G_2(929)$, $L_{2}(5503^{2})$, $S_{4}(5503)$, $U_{3}(5503)$, $U_{4}(5503)$\\
\SetRow{lightgray!40} 5507 & 17 & $L_{2}(5507^{2})$, $S_{4}(5507)$, $U_{3}(5507)$, $U_{4}(5507)$\\
5527 & 11 & $L_{2}(877^{3})$, $G_2(877)$, $U_{3}(877)$, $L_{2}(4651^{3})$, $G_2(4651)$, $U_{3}(4651)$\\
\SetRow{lightgray!40} 5531 & 30 & $O^-_{10}(239)$, $U_{5}(239)$, $U_{6}(239)$\\
5557 & 14 & $L_{2}(3079^{2})$, $S_{4}(3079)$, $U_{4}(3079)$, $L_{2}(5557^{2})$, $S_{4}(5557)$, $U_{3}(5557)$, $U_{4}(5557)$\\
\SetRow{lightgray!40} 5569 & 8 & $L_{2}(2243^{3})$, $G_2(2243)$, $U_{3}(2243)$\\
5573 & 12 & $L_{4}(2017)$, $L_{2}(2017^{2})$, $S_{4}(2017)$\\
\SetRow{lightgray!40} 5581 & 19 & $L_{5}(53)$, $L_{6}(53)$, $L_{2}(53^{5})$, $U_{3}(2459)$, $L_{3}(5581)$, $L_{4}(5581)$, $L_{2}(5581^{2})$, $S_{4}(5581)$\\
5591 & 34 & $U_{3}(5591)$\\
\SetRow{lightgray!40} 5623 & 19 & $L_{2}(5623^{2})$, $S_{4}(5623)$\\
5639 & 4 & $L_{3}(5639)$\\
\SetRow{lightgray!40} 5641 & 17 & $L_{4}(1429)$, $L_{2}(1429^{2})$, $L_{3}(1429^{2})$, $S_{4}(1429)$, $S_{6}(1429)$, $O_{7}(1429)$, $O^+_{8}(1429)$, $U_{4}(1429)$, $L_{2}(5641^{2})$, $S_{4}(5641)$\\
5647 & 14 & $L_{3}(853)$, $L_{4}(853)$, $L_{3}(853^{2})$, $L_{2}(853^{3})$, $S_{6}(853)$, $O_{7}(853)$, $O^+_{8}(853)$, $G_2(853)$, $L_{3}(4793)$\\
\SetRow{lightgray!40} 5653 & 6 & $U_{3}(17^{3})$\\
5657 & 4 & $U_{3}(5657)$\\
\SetRow{lightgray!40} 5659 & 15 & $L_{3}(5659)$, $L_{2}(5659^{3})$, $G_2(5659)$, $U_{3}(5659)$\\
5683 & 8 & $U_{3}(5683)$\\
\SetRow{lightgray!40} 5693 & 13 & $L_{2}(1193^{2})$, $S_{4}(1193)$, $U_{4}(1193)$, $U_{3}(5693)$\\
5701 & 13 & $L_{2}(5701^{2})$, $S_{4}(5701)$\\
\SetRow{lightgray!40} 5711 & 8 & $L_{3}(5711)$\\
5717 & 23 & $L_{2}(3301^{2})$, $S_{4}(3301)$\\
\SetRow{lightgray!40} 5743 & 11 & $L_{3}(5743)$, $L_{4}(5743)$, $L_{2}(5743^{2})$, $S_{4}(5743)$\\
5749 & 44 & $L_{3}(331^{2})$, $L_{2}(331^{3})$, $S_{6}(331)$, $O_{7}(331)$, $O^+_{8}(331)$, $G_2(331)$, $U_{3}(331)$, $U_{4}(331)$, $L_{2}(4943^{2})$, $S_{4}(4943)$, $L_{2}(5419^{3})$, $G_2(5419)$, $U_{3}(5419)$\\
\SetRow{lightgray!40} 5779 & 8 & $L_{3}(2851)$, $L_{3}(2927)$, $L_{3}(5779)$\\
5791 & 20 & $L_{3}(4219^{2})$, $L_{2}(4219^{3})$, $S_{6}(4219)$, $O_{7}(4219)$, $O^+_{8}(4219)$, $G_2(4219)$, $U_{3}(4219)$, $U_{4}(4219)$, $L_{3}(5791)$\\
\SetRow{lightgray!40} 5801 & 8 & $L_{3}(5801)$\\
5807 & 11 & $L_{3}(5807)$, $L_{4}(5807)$, $L_{2}(5807^{2})$, $S_{4}(5807)$\\
\SetRow{lightgray!40} 5813 & 11 & $L_{2}(5813^{2})$, $S_{4}(5813)$\\
5821 & 11 & $L_{2}(3673^{3})$, $G_2(3673)$, $U_{3}(3673)$, $L_{3}(5821)$\\
\SetRow{lightgray!40} 5827 & 17 & $L_{3}(5827)$, $L_{2}(5827^{3})$, $G_2(5827)$, $U_{3}(5827)$\\
5843 & 9 & $L_{2}(5843^{2})$, $S_{4}(5843)$\\
\SetRow{lightgray!40} 5849 & 4 & $U_{3}(5849)$\\
5851 & 16 & $L_{3}(577)$, $L_{4}(577)$, $L_{3}(577^{2})$, $L_{2}(577^{3})$, $S_{6}(577)$, $O_{7}(577)$, $O^+_{8}(577)$, $G_2(577)$, $L_{3}(5273)$\\
\SetRow{lightgray!40} 5857 & 8 & $L_{2}(4547^{2})$, $S_{4}(4547)$, $U_{3}(5857)$\\
5861 & 10 & $L_{2}(5107^{2})$, $S_{4}(5107)$, $L_{3}(5861)$\\
\SetRow{lightgray!40} 5867 & 5 & $L_{2}(5867^{2})$, $S_{4}(5867)$\\
5869 & 15 & $L_{3}(5869)$, $L_{4}(5869)$, $L_{2}(5869^{2})$, $S_{4}(5869)$\\
\SetRow{lightgray!40} 5879 & 4 & $U_{3}(5879)$\\
5881 & 28 & $L_{2}(277^{3})$, $G_2(277)$, $U_{3}(277)$, $L_{4}(4783)$, $L_{2}(4783^{2})$, $L_{3}(4783^{2})$, $S_{4}(4783)$, $S_{6}(4783)$, $O_{7}(4783)$, $O^+_{8}(4783)$, $U_{4}(4783)$\\
\SetRow{lightgray!40} 5923 & 7 & $L_{2}(5923^{2})$, $S_{4}(5923)$\\
5927 & 14 & $U_{3}(5927)$\\
\SetRow{lightgray!40} 6007 & 7 & $L_{2}(6007^{2})$, $S_{4}(6007)$\\
6029 & 12 & $L_{2}(1801^{2})$, $S_{4}(1801)$, $U_{3}(6029)$\\
\SetRow{lightgray!40} 6037 & 16 & $L_{3}(509)$, $L_{4}(509)$, $L_{3}(509^{2})$, $L_{2}(509^{3})$, $S_{6}(509)$, $O_{7}(509)$, $O^+_{8}(509)$, $G_2(509)$, $L_{3}(5527)$\\
6043 & 9 & $L_{3}(4327)$, $L_{4}(4327)$, $L_{2}(6043^{2})$, $S_{4}(6043)$\\
\SetRow{lightgray!40} 6053 & 18 & $L_{2}(3221^{2})$, $S_{4}(3221)$, $U_{4}(3221)$\\
6067 & 11 & $L_{3}(6067)$, $L_{4}(6067)$, $L_{2}(6067^{2})$, $S_{4}(6067)$\\
\SetRow{lightgray!40} 6073 & 8 & $U_{3}(4231)$\\
6079 & 20 & $L_{3}(1553)$, $L_{4}(1553)$, $L_{3}(1553^{2})$, $L_{2}(1553^{3})$, $S_{6}(1553)$, $O_{7}(1553)$, $O^+_{8}(1553)$, $G_2(1553)$, $L_{3}(6079)$\\
\SetRow{lightgray!40} 6089 & 5 & $L_{2}(6089^{2})$, $S_{4}(6089)$\\
6091 & 14 & $U_{3}(5347)$, $U_{4}(5347)$, $L_{3}(6091)$\\
\SetRow{lightgray!40} 6121 & 27 & $L_{3}(1153^{2})$, $L_{2}(1153^{3})$, $S_{6}(1153)$, $O_{7}(1153)$, $O^+_{8}(1153)$, $G_2(1153)$, $U_{3}(1153)$, $U_{4}(1153)$, $L_{3}(4969^{2})$, $L_{2}(4969^{3})$, $S_{6}(4969)$, $O_{7}(4969)$, $O^+_{8}(4969)$, $G_2(4969)$, $U_{3}(4969)$, $U_{4}(4969)$\\
6133 & 15 & $L_{3}(6133)$, $L_{4}(6133)$, $L_{2}(6133^{2})$, $S_{4}(6133)$\\
\SetRow{lightgray!40} 6143 & 11 & $L_{2}(6143^{2})$, $S_{4}(6143)$\\
6151 & 15 & $L_{2}(6151^{2})$, $S_{4}(6151)$\\
\SetRow{lightgray!40} 6163 & 28 & $L_{3}(79^{2})$, $L_{2}(79^{3})$, $L_{4}(79^{2})$, $S_{6}(79)$, $S_{8}(79)$, $O_{7}(79)$, $O_{9}(79)$, $O^+_{8}(79)$, $O^-_{8}(79)$, $O^-_{10}(79)$, $G_2(79)$, $U_{3}(79)$, $U_{4}(79)$, $U_{5}(79)$, $U_{6}(79)$, $L_{2}(6163^{2})$, $S_{4}(6163)$\\
6173 & 29 & $L_{2}(2447^{2})$, $S_{4}(2447)$, $L_{2}(6173^{2})$, $S_{4}(6173)$\\
\SetRow{lightgray!40} 6197 & 5 & $L_{2}(6197^{2})$, $S_{4}(6197)$\\
6211 & 22 & $L_{3}(137^{2})$, $L_{2}(137^{3})$, $L_{4}(137^{2})$, $S_{6}(137)$, $S_{8}(137)$, $O_{7}(137)$, $O_{9}(137)$, $O^+_{8}(137)$, $O^-_{8}(137)$, $O^-_{10}(137)$, $G_2(137)$, $U_{3}(137)$, $U_{4}(137)$, $U_{5}(137)$, $U_{6}(137)$\\
\SetRow{lightgray!40} 6217 & 10 & $L_{3}(2459)$, $L_{2}(2459^{3})$, $G_2(2459)$, $L_{2}(6217^{2})$, $S_{4}(6217)$\\
6229 & 23 & $L_{4}(1451)$, $L_{2}(1451^{2})$, $S_{4}(1451)$, $L_{3}(6229)$\\
\SetRow{lightgray!40} 6247 & 14 & $U_{3}(3931)$, $L_{2}(6247^{2})$, $S_{4}(6247)$\\
6257 & 11 & $L_{2}(4673^{2})$, $S_{4}(4673)$, $L_{2}(6257^{2})$, $S_{4}(6257)$\\
\SetRow{lightgray!40} 6263 & 9 & $L_{2}(6263^{2})$, $S_{4}(6263)$\\
6269 & 9 & $L_{2}(1523^{2})$, $S_{4}(1523)$, $L_{3}(6269)$, $L_{4}(6269)$, $L_{2}(6269^{2})$, $S_{4}(6269)$\\
\SetRow{lightgray!40} 6277 & 25 & $L_{4}(1033)$, $L_{2}(1033^{2})$, $L_{3}(1033^{2})$, $S_{4}(1033)$, $S_{6}(1033)$, $O_{7}(1033)$, $O^+_{8}(1033)$, $U_{4}(1033)$, $U_{3}(2309)$, $U_{4}(2309)$, $L_{3}(6277)$, $L_{4}(6277)$, $L_{2}(6277^{2})$, $S_{4}(6277)$\\
6287 & 17 & $L_{3}(6287)$, $L_{2}(6287^{3})$, $G_2(6287)$, $U_{3}(6287)$\\
\SetRow{lightgray!40} 6299 & 4 & $U_{3}(6299)$\\
6301 & 19 & $L_{3}(3323)$, $L_{4}(3323)$, $L_{3}(3323^{2})$, $L_{2}(3323^{3})$, $S_{6}(3323)$, $O_{7}(3323)$, $O^+_{8}(3323)$, $G_2(3323)$\\
\SetRow{lightgray!40} 6311 & 9 & $L_{2}(6311^{2})$, $S_{4}(6311)$\\
6317 & 8 & $L_{3}(6317)$\\
\SetRow{lightgray!40} 6323 & 11 & $L_{3}(6323)$, $L_{4}(6323)$, $L_{2}(6323^{2})$, $S_{4}(6323)$\\
6337 & 9 & $L_{2}(6337^{2})$, $S_{4}(6337)$\\
\SetRow{lightgray!40} 6343 & 13 & $L_{3}(557)$, $L_{4}(557)$\\
6353 & 11 & $L_{3}(6353)$, $L_{4}(6353)$, $L_{2}(6353^{2})$, $S_{4}(6353)$\\
\SetRow{lightgray!40} 6359 & 4 & $U_{3}(6359)$\\
6367 & 8 & $L_{3}(769)$\\
\SetRow{lightgray!40} 6373 & 11 & $L_{2}(1879^{2})$, $S_{4}(1879)$, $U_{4}(1879)$, $L_{3}(5749)$\\
6379 & 19 & $L_{3}(3373^{2})$, $L_{2}(3373^{3})$, $S_{6}(3373)$, $O_{7}(3373)$, $O^+_{8}(3373)$, $G_2(3373)$, $U_{3}(3373)$, $U_{4}(3373)$\\
\SetRow{lightgray!40} 6389 & 12 & $L_{2}(4297^{2})$, $S_{4}(4297)$, $U_{4}(4297)$\\
6397 & 27 & $L_{2}(6397^{2})$, $S_{4}(6397)$\\
\SetRow{lightgray!40} 6421 & 8 & $L_{3}(6421)$\\
6427 & 24 & $L_{3}(6427)$\\
\SetRow{lightgray!40} 6449 & 4 & $L_{3}(6449)$\\
6451 & 21 & $L_{2}(6451^{2})$, $S_{4}(6451)$\\
\SetRow{lightgray!40} 6469 & 8 & $U_{3}(4993)$, $L_{2}(6469^{2})$, $S_{4}(6469)$\\
6473 & 10 & $U_{3}(6473)$\\
\SetRow{lightgray!40} 6481 & 50 & $L_{12}(9)$, $L_{6}(3^{4})$, $L_{4}(3^{6})$, $L_{3}(3^{8})$, $L_{2}(3^{12})$, $L_{5}(3^{6})$, $L_{2}(3^{24})$, $S_{24}(3)$, $S_{12}(9)$, $S_{8}(27)$, $S_{6}(3^{4})$, $S_{4}(3^{6})$, $S_{10}(27)$, $S_{4}(3^{12})$, $O_{25}(3)$, $O_{13}(9)$, $O_{9}(27)$, $O_{7}(3^{4})$, $O_{11}(27)$, $O^+_{14}(9)$, $O^+_{10}(27)$, $O^+_{8}(3^{4})$, $O^+_{12}(27)$, $O^-_{8}(27)$, $O^-_{12}(9)$, $O^-_{24}(3)$, $O^-_{10}(27)$, $G_2(3^{4})$, $G_2(3^{8})$, $F_4(9)$, $E_6(9)$, $E_8(3)$, $^3D_4(9)$, $^3D_4(3^{4})$, $U_{3}(3^{4})$, $U_{4}(3^{4})$, $U_{8}(27)$, $L_{2}(6481^{2})$, $S_{4}(6481)$\\
6491 & 32 & $U_{3}(6491)$\\
\SetRow{lightgray!40} 6521 & 13 & $L_{2}(4157^{2})$, $S_{4}(4157)$, $L_{2}(6521^{2})$, $S_{4}(6521)$\\
6529 & 32 & $L_{3}(491)$, $L_{4}(491)$, $L_{3}(491^{2})$, $L_{2}(491^{3})$, $S_{6}(491)$, $O_{7}(491)$, $O^+_{8}(491)$, $G_2(491)$, $L_{4}(2311)$, $L_{2}(2311^{2})$, $S_{4}(2311)$, $L_{3}(6037)$, $U_{3}(6529)$\\
\SetRow{lightgray!40} 6547 & 8 & $U_{3}(2333)$, $U_{4}(2333)$, $U_{3}(6547)$\\
6551 & 5 & $L_{2}(6551^{2})$, $S_{4}(6551)$\\
\SetRow{lightgray!40} 6553 & 13 & $L_{2}(6553^{2})$, $S_{4}(6553)$\\
6569 & 5 & $L_{2}(6569^{2})$, $S_{4}(6569)$\\
\SetRow{lightgray!40} 6571 & 9 & $L_{2}(6571^{2})$, $S_{4}(6571)$\\
6577 & 13 & $L_{3}(353)$, $L_{4}(353)$, $L_{3}(353^{2})$, $L_{2}(353^{3})$, $S_{6}(353)$, $O_{7}(353)$, $O^+_{8}(353)$, $G_2(353)$\\
\SetRow{lightgray!40} 6581 & 20 & $L_{3}(6581)$\\
6607 & 14 & $L_{3}(6607)$\\
\SetRow{lightgray!40} 6619 & 20 & $L_{3}(569)$\\
6637 & 18 & $L_{3}(6637)$\\
\SetRow{lightgray!40} 6659 & 4 & $L_{3}(6659)$\\
6661 & 15 & $L_{2}(6661^{2})$, $S_{4}(6661)$\\
\SetRow{lightgray!40} 6673 & 10 & $L_{2}(2437^{2})$, $S_{4}(2437)$, $L_{3}(5279)$\\
6679 & 12 & $U_{3}(5737)$\\
\SetRow{lightgray!40} 6691 & 13 & $L_{2}(6691^{2})$, $S_{4}(6691)$\\
6701 & 7 & $L_{4}(1721)$, $L_{2}(1721^{2})$, $S_{4}(1721)$, $L_{3}(6701)$\\
\SetRow{lightgray!40} 6703 & 10 & $U_{3}(1481)$, $U_{4}(1481)$, $L_{3}(6703)$\\
6733 & 13 & $L_{3}(619)$, $L_{2}(619^{3})$, $G_2(619)$, $L_{3}(6113)$, $L_{3}(6733)$, $L_{2}(6733^{3})$, $G_2(6733)$, $U_{3}(6733)$\\
\SetRow{lightgray!40} 6737 & 34 & $L_{4}(2393)$, $L_{2}(2393^{2})$, $L_{3}(2393^{2})$, $S_{4}(2393)$, $S_{6}(2393)$, $O_{7}(2393)$, $O^+_{8}(2393)$, $U_{4}(2393)$, $L_{3}(6737)$\\
6763 & 18 & $L_{3}(6763)$\\
\SetRow{lightgray!40} 6779 & 4 & $U_{3}(6779)$\\
6781 & 14 & $L_{2}(2927^{3})$, $G_2(2927)$, $U_{3}(2927)$\\
\SetRow{lightgray!40} 6793 & 15 & $L_{2}(709^{2})$, $S_{4}(709)$, $U_{4}(709)$, $L_{3}(6793)$\\
6803 & 25 & $L_{3}(6803)$, $L_{2}(6803^{3})$, $G_2(6803)$, $U_{3}(6803)$\\
\SetRow{lightgray!40} 6823 & 6 & $U_{3}(6823)$\\
6827 & 5 & $L_{2}(6827^{2})$, $S_{4}(6827)$\\
\SetRow{lightgray!40} 6829 & 8 & $L_{2}(5233^{2})$, $S_{4}(5233)$, $L_{3}(6829)$\\
6833 & 11 & $L_{2}(1307^{2})$, $S_{4}(1307)$\\
\SetRow{lightgray!40} 6841 & 23 & $L_{2}(53^{6})$, $S_{4}(53^{3})$, $G_2(53^{2})$, $^3D_4(53)$, $U_{3}(53^{2})$, $L_{3}(6841)$\\
6863 & 19 & $L_{3}(6863)$, $L_{4}(6863)$, $L_{2}(6863^{2})$, $L_{3}(6863^{2})$, $L_{2}(6863^{3})$, $S_{4}(6863)$, $S_{6}(6863)$, $O_{7}(6863)$, $O^+_{8}(6863)$, $G_2(6863)$, $U_{3}(6863)$, $U_{4}(6863)$\\
\SetRow{lightgray!40} 6871 & 14 & $L_{3}(6871)$\\
6899 & 10 & $L_{3}(6899)$\\
\SetRow{lightgray!40} 6907 & 12 & $L_{2}(5051^{3})$, $G_2(5051)$, $U_{3}(5051)$, $L_{3}(6907)$, $L_{2}(6907^{3})$, $G_2(6907)$, $U_{3}(6907)$\\
6911 & 9 & $L_{2}(6911^{2})$, $S_{4}(6911)$\\
\SetRow{lightgray!40} 6917 & 41 & $L_{4}(263)$, $L_{2}(263^{2})$, $L_{3}(263^{2})$, $S_{4}(263)$, $S_{6}(263)$, $O_{7}(263)$, $O^+_{8}(263)$, $U_{4}(263)$, $L_{2}(6917^{2})$, $S_{4}(6917)$\\
6949 & 12 & $U_{3}(6949)$\\
\SetRow{lightgray!40} 6959 & 5 & $L_{2}(6959^{2})$, $S_{4}(6959)$\\
6961 & 9 & $U_{3}(727^{2})$, $L_{3}(6961)$\\
\SetRow{lightgray!40} 6967 & 7 & $U_{3}(383)$, $U_{3}(6967)$\\
6971 & 8 & $L_{3}(6971)$\\
\SetRow{lightgray!40} 6977 & 11 & $L_{2}(2063^{2})$, $S_{4}(2063)$, $L_{2}(6977^{2})$, $S_{4}(6977)$\\
6983 & 13 & $L_{2}(6983^{2})$, $S_{4}(6983)$, $U_{3}(6983)$, $U_{4}(6983)$\\
\SetRow{lightgray!40} 6991 & 8 & $L_{3}(1381)$\\
6997 & 7 & $U_{3}(2909)$, $L_{3}(6997)$\\
\SetRow{lightgray!40} 7027 & 22 & $L_{3}(523)$, $L_{4}(523)$, $L_{3}(523^{2})$, $L_{2}(523^{3})$, $S_{6}(523)$, $O_{7}(523)$, $O^+_{8}(523)$, $G_2(523)$, $U_{6}(523)$\\
7039 & 8 & $L_{2}(6737^{3})$, $G_2(6737)$, $U_{3}(6737)$\\
\SetRow{lightgray!40} 7057 & 16 & $L_{3}(6911)$, $L_{4}(6911)$, $U_{3}(7057)$\\
7079 & 27 & $L_{2}(7079^{2})$, $S_{4}(7079)$\\
\SetRow{lightgray!40} 7103 & 12 & $L_{2}(7103^{2})$, $S_{4}(7103)$, $U_{3}(7103)$, $U_{4}(7103)$, $U_{5}(7103)$\\
7121 & 10 & $L_{2}(6343^{2})$, $S_{4}(6343)$, $L_{3}(7121)$\\
\SetRow{lightgray!40} 7127 & 5 & $L_{2}(7127^{2})$, $S_{4}(7127)$\\
7129 & 27 & $L_{3}(1249)$, $L_{3}(5879)$, $L_{2}(5879^{3})$, $G_2(5879)$\\
\SetRow{lightgray!40} 7151 & 10 & $L_{3}(7151)$\\
7159 & 20 & $L_{3}(7159)$\\
\SetRow{lightgray!40} 7177 & 20 & $L_{3}(2039^{2})$, $L_{2}(2039^{3})$, $S_{6}(2039)$, $O_{7}(2039)$, $O^+_{8}(2039)$, $G_2(2039)$, $U_{3}(2039)$, $U_{4}(2039)$, $L_{3}(7177)$\\
7187 & 11 & $L_{3}(7187)$, $L_{4}(7187)$, $L_{2}(7187^{2})$, $S_{4}(7187)$\\
\SetRow{lightgray!40} 7193 & 22 & $L_{4}(967)$, $L_{2}(967^{2})$, $S_{4}(967)$, $L_{3}(7193)$, $L_{4}(7193)$, $L_{2}(7193^{2})$, $S_{4}(7193)$\\
7207 & 9 & $L_{2}(7207^{2})$, $S_{4}(7207)$, $U_{3}(7207)$, $U_{4}(7207)$\\
\SetRow{lightgray!40} 7213 & 16 & $L_{4}(1999)$, $L_{2}(1999^{2})$, $L_{3}(1999^{2})$, $S_{4}(1999)$, $S_{6}(1999)$, $O_{7}(1999)$, $O^+_{8}(1999)$, $U_{4}(1999)$, $U_{3}(7213)$\\
7219 & 13 & $L_{3}(4493)$, $L_{4}(4493)$\\
\SetRow{lightgray!40} 7237 & 11 & $L_{2}(1831^{3})$, $G_2(1831)$, $U_{3}(1831)$, $U_{3}(5407)$\\
7253 & 33 & $L_{2}(7253^{2})$, $S_{4}(7253)$\\
\SetRow{lightgray!40} 7283 & 16 & $L_{3}(7283)$\\
7297 & 15 & $L_{3}(3761)$, $L_{4}(3761)$, $L_{2}(7297^{2})$, $S_{4}(7297)$\\
\SetRow{lightgray!40} 7307 & 5 & $L_{2}(7307^{2})$, $S_{4}(7307)$\\
7309 & 14 & $U_{3}(7309)$\\
\SetRow{lightgray!40} 7321 & 25 & $L_{4}(11^{2})$, $L_{2}(11^{4})$, $L_{3}(11^{4})$, $S_{8}(11)$, $S_{4}(11^{2})$, $S_{6}(11^{2})$, $O_{9}(11)$, $O_{7}(11^{2})$, $O^+_{10}(11)$, $O^+_{8}(11^{2})$, $O^-_{8}(11)$, $F_4(11)$, $U_{4}(11^{2})$, $U_{3}(7013)$\\
7333 & 23 & $L_{2}(2909^{2})$, $S_{4}(2909)$, $U_{4}(2909)$, $L_{2}(4271^{3})$, $G_2(4271)$, $U_{3}(4271)$\\
\SetRow{lightgray!40} 7351 & 28 & $L_{6}(149)$, $L_{3}(149^{2})$, $L_{2}(149^{3})$, $S_{6}(149)$, $O_{7}(149)$, $O^+_{8}(149)$, $G_2(149)$, $U_{3}(149)$, $U_{4}(149)$\\
7369 & 33 & $L_{2}(607^{2})$, $S_{4}(607)$, $L_{4}(3373^{2})$, $L_{2}(3373^{4})$, $S_{8}(3373)$, $S_{4}(3373^{2})$, $O_{9}(3373)$, $O^-_{8}(3373)$\\
\SetRow{lightgray!40} 7433 & 22 & $L_{2}(983^{2})$, $S_{4}(983)$, $U_{4}(983)$\\
7459 & 27 & $L_{3}(229^{2})$, $L_{2}(229^{3})$, $S_{6}(229)$, $O_{7}(229)$, $O^+_{8}(229)$, $G_2(229)$, $U_{3}(229)$, $U_{4}(229)$\\
\SetRow{lightgray!40} 7477 & 7 & $U_{3}(3469)$, $U_{4}(3469)$\\
7481 & 9 & $L_{2}(6073^{2})$, $S_{4}(6073)$\\
\SetRow{lightgray!40} 7489 & 13 & $L_{3}(2467)$, $L_{3}(5021)$\\
7499 & 10 & $L_{3}(7499)$\\
\SetRow{lightgray!40} 7507 & 14 & $U_{3}(607)$, $U_{4}(607)$, $U_{3}(7507)$\\
7529 & 13 & $L_{3}(7529)$, $L_{2}(7529^{3})$, $G_2(7529)$, $U_{3}(7529)$\\
\SetRow{lightgray!40} 7537 & 9 & $L_{4}(1049)$, $L_{2}(1049^{2})$, $S_{4}(1049)$, $L_{3}(7537)$\\
7549 & 21 & $L_{2}(23^{6})$, $S_{4}(23^{3})$, $G_2(23^{2})$, $^3D_4(23)$, $U_{3}(23^{2})$, $U_{4}(23^{3})$, $L_{2}(7549^{2})$, $S_{4}(7549)$, $U_{3}(7549)$, $U_{4}(7549)$\\
\SetRow{lightgray!40} 7559 & 4 & $L_{3}(7559)$\\
7561 & 15 & $U_{3}(6263)$, $U_{4}(6263)$\\
\SetRow{lightgray!40} 7577 & 11 & $L_{4}(6037)$, $L_{2}(6037^{2})$, $S_{4}(6037)$, $L_{3}(7577)$\\
7583 & 8 & $U_{3}(7583)$\\
\SetRow{lightgray!40} 7603 & 9 & $L_{2}(7603^{2})$, $S_{4}(7603)$, $U_{3}(7603)$, $U_{4}(7603)$\\
7607 & 19 & $L_{2}(7607^{2})$, $S_{4}(7607)$, $U_{3}(7607)$, $U_{4}(7607)$\\
\SetRow{lightgray!40} 7621 & 26 & $L_{2}(5^{15})$, $G_2(5^{5})$, $U_{5}(5^{3})$, $U_{3}(5^{5})$, $U_{6}(5^{3})$, $U_{7}(5^{3})$, $L_{3}(7621)$\\
7639 & 10 & $L_{3}(4663)$, $L_{3}(7639)$, $L_{4}(7639)$, $L_{2}(7639^{2})$, $S_{4}(7639)$\\
\SetRow{lightgray!40} 7649 & 23 & $L_{2}(7649^{2})$, $S_{4}(7649)$\\
7669 & 9 & $L_{2}(7669^{2})$, $S_{4}(7669)$, $U_{3}(7669)$, $U_{4}(7669)$\\
\SetRow{lightgray!40} 7673 & 19 & $L_{4}(277)$, $L_{2}(277^{2})$, $L_{3}(277^{2})$, $S_{4}(277)$, $S_{6}(277)$, $O_{7}(277)$, $O^+_{8}(277)$, $U_{4}(277)$, $U_{5}(277)$, $U_{6}(277)$\\
7681 & 11 & $L_{2}(6997^{3})$, $G_2(6997)$, $U_{3}(6997)$, $U_{3}(7681)$\\
\SetRow{lightgray!40} 7687 & 8 & $U_{3}(5413)$, $U_{4}(5413)$, $L_{3}(7687)$\\
7699 & 10 & $L_{3}(2269)$, $L_{2}(7699^{2})$, $S_{4}(7699)$, $U_{3}(7699)$, $U_{4}(7699)$\\
\SetRow{lightgray!40} 7703 & 16 & $L_{3}(7703)$\\
7717 & 10 & $L_{2}(2953^{2})$, $S_{4}(2953)$, $U_{4}(2953)$\\
\SetRow{lightgray!40} 7723 & 6 & $U_{3}(7723)$\\
7727 & 16 & $L_{3}(7727)$\\
\SetRow{lightgray!40} 7753 & 6 & $L_{3}(7349)$\\
7759 & 35 & $L_{3}(1759)$, $L_{2}(1759^{3})$, $G_2(1759)$, $U_{3}(7759)$\\
\SetRow{lightgray!40} 7789 & 15 & $L_{3}(233)$, $L_{4}(233)$, $L_{3}(233^{2})$, $L_{2}(233^{3})$, $S_{6}(233)$, $O_{7}(233)$, $O^+_{8}(233)$, $G_2(233)$, $L_{2}(7789^{2})$, $S_{4}(7789)$\\
7841 & 17 & $L_{2}(7643^{2})$, $S_{4}(7643)$, $L_{2}(7841^{2})$, $S_{4}(7841)$\\
\SetRow{lightgray!40} 7853 & 23 & $L_{4}(1759)$, $L_{2}(1759^{2})$, $L_{3}(1759^{2})$, $S_{4}(1759)$, $S_{6}(1759)$, $O_{7}(1759)$, $O^+_{8}(1759)$, $U_{4}(1759)$\\
7867 & 9 & $L_{2}(7867^{2})$, $S_{4}(7867)$\\
\SetRow{lightgray!40} 7873 & 8 & $L_{4}(4283)$, $L_{2}(4283^{2})$, $S_{4}(4283)$\\
7877 & 4 & $U_{3}(7877)$\\
\SetRow{lightgray!40} 7879 & 9 & $L_{2}(1367^{3})$, $G_2(1367)$, $U_{3}(1367)$, $L_{3}(7879)$\\
7883 & 23 & $L_{2}(7883^{2})$, $S_{4}(7883)$, $U_{3}(7883)$, $U_{4}(7883)$\\
\SetRow{lightgray!40} 7907 & 14 & $L_{3}(7907)$\\
7919 & 10 & $U_{3}(7919)$\\
\SetRow{lightgray!40} 7927 & 9 & $L_{2}(7927^{2})$, $S_{4}(7927)$\\
7933 & 8 & $L_{3}(5927)$, $L_{2}(5927^{3})$, $G_2(5927)$\\
\SetRow{lightgray!40} 7949 & 4 & $U_{3}(7949)$\\
7963 & 33 & $L_{2}(7963^{2})$, $S_{4}(7963)$\\
\SetRow{lightgray!40} 7993 & 18 & $U_{3}(7993)$\\
8009 & 11 & $L_{4}(283)$, $L_{2}(283^{2})$, $L_{3}(283^{2})$, $S_{4}(283)$, $S_{6}(283)$, $O_{7}(283)$, $O^+_{8}(283)$, $U_{4}(283)$\\
\SetRow{lightgray!40} 8011 & 22 & $L_{3}(89)$, $L_{4}(89)$, $L_{5}(89)$, $L_{6}(89)$, $L_{3}(89^{2})$, $L_{2}(89^{3})$, $L_{2}(89^{6})$, $S_{6}(89)$, $S_{4}(89^{3})$, $O_{7}(89)$, $O^+_{8}(89)$, $G_2(89)$, $G_2(89^{2})$, $^3D_4(89)$, $U_{3}(8011)$\\
8039 & 27 & $L_{3}(8039)$, $L_{4}(8039)$, $L_{2}(8039^{2})$, $L_{3}(8039^{2})$, $L_{2}(8039^{3})$, $S_{4}(8039)$, $S_{6}(8039)$, $O_{7}(8039)$, $O^+_{8}(8039)$, $G_2(8039)$, $U_{3}(8039)$, $U_{4}(8039)$\\
\SetRow{lightgray!40} 8069 & 14 & $U_{3}(8069)$\\
8089 & 8 & $L_{4}(2293)$, $L_{2}(2293^{2})$, $S_{4}(2293)$\\
\SetRow{lightgray!40} 8093 & 13 & $L_{2}(8093^{2})$, $S_{4}(8093)$, $U_{3}(8093)$, $U_{4}(8093)$\\
8101 & 15 & $L_{2}(8011^{2})$, $S_{4}(8011)$, $U_{4}(8011)$, $L_{3}(8101)$\\
\SetRow{lightgray!40} 8111 & 8 & $U_{3}(8111)$\\
8117 & 11 & $L_{2}(1733^{2})$, $S_{4}(1733)$, $L_{2}(8117^{2})$, $S_{4}(8117)$\\
\SetRow{lightgray!40} 8147 & 16 & $L_{3}(8147)$\\
8161 & 13 & $L_{3}(2903)$, $L_{4}(2903)$, $L_{2}(8161^{2})$, $S_{4}(8161)$, $U_{3}(8161)$, $U_{4}(8161)$\\
\SetRow{lightgray!40} 8179 & 15 & $U_{3}(1097)$, $U_{3}(8179)$\\
8191 & 44 & $L_{13}(2)$, $L_{14}(2)$, $L_{15}(2)$, $L_{16}(2)$, $L_{13}(4)$, $L_{2}(2^{13})$, $L_{14}(4)$, $L_{15}(4)$, $L_{2}(2^{26})$, $S_{26}(2)$, $S_{28}(2)$, $S_{30}(2)$, $S_{4}(2^{13})$, $O^+_{26}(2)$, $O^+_{28}(2)$, $O^+_{30}(2)$, $O^+_{32}(2)$, $O^-_{28}(2)$, $O^-_{30}(2)$, $L_{2}(8101^{3})$, $G_2(8101)$, $U_{3}(8101)$, $L_{2}(8191^{2})$, $S_{4}(8191)$, $Sz(2^{13})$\\
\SetRow{lightgray!40} 8209 & 17 & $L_{3}(3^{9})$, $L_{2}(2383^{4})$, $S_{4}(2383^{2})$, $L_{3}(4943)$, $L_{4}(4943)$, $L_{3}(8209)$\\
8219 & 4 & $L_{3}(8219)$\\
\SetRow{lightgray!40} 8221 & 15 & $L_{2}(8221^{2})$, $S_{4}(8221)$, $U_{3}(8221)$, $U_{4}(8221)$\\
8237 & 9 & $L_{2}(8237^{2})$, $S_{4}(8237)$\\
\SetRow{lightgray!40} 8243 & 22 & $L_{3}(8243)$\\
8263 & 9 & $U_{3}(241)$, $U_{4}(241)$\\
\SetRow{lightgray!40} 8269 & 19 & $L_{3}(157)$, $L_{4}(157)$, $L_{3}(157^{2})$, $L_{2}(157^{3})$, $S_{6}(157)$, $O_{7}(157)$, $O^+_{8}(157)$, $G_2(157)$, $L_{4}(643)$, $L_{2}(643^{2})$, $S_{4}(643)$, $L_{3}(8111)$, $L_{2}(8111^{3})$, $G_2(8111)$\\
8273 & 16 & $U_{3}(8273)$\\
\SetRow{lightgray!40} 8287 & 9 & $L_{2}(569^{3})$, $G_2(569)$, $U_{3}(569)$, $U_{3}(8287)$\\
8291 & 7 & $L_{3}(8291)$, $L_{2}(8291^{3})$, $G_2(8291)$, $U_{3}(8291)$\\
\SetRow{lightgray!40} 8297 & 27 & $L_{3}(8297)$, $L_{4}(8297)$, $L_{2}(8297^{2})$, $L_{3}(8297^{2})$, $L_{2}(8297^{3})$, $S_{4}(8297)$, $S_{6}(8297)$, $O_{7}(8297)$, $O^+_{8}(8297)$, $G_2(8297)$, $U_{3}(8297)$, $U_{4}(8297)$\\
8317 & 15 & $L_{2}(8317^{2})$, $S_{4}(8317)$\\
\SetRow{lightgray!40} 8329 & 26 & $L_{3}(8329)$\\
8363 & 11 & $L_{2}(8363^{2})$, $S_{4}(8363)$, $U_{3}(8363)$, $U_{4}(8363)$\\
\SetRow{lightgray!40} 8369 & 12 & $L_{4}(7703)$, $L_{2}(7703^{2})$, $S_{4}(7703)$\\
8377 & 13 & $L_{2}(8377^{2})$, $S_{4}(8377)$\\
\SetRow{lightgray!40} 8387 & 5 & $L_{2}(8387^{2})$, $S_{4}(8387)$\\
8389 & 37 & $L_{3}(691)$, $L_{4}(691)$, $L_{2}(3449^{2})$, $S_{4}(3449)$, $U_{4}(3449)$, $L_{3}(8389)$\\
\SetRow{lightgray!40} 8423 & 9 & $L_{2}(8423^{2})$, $S_{4}(8423)$\\
8429 & 7 & $L_{3}(8429)$, $L_{2}(8429^{3})$, $G_2(8429)$, $U_{3}(8429)$\\
\SetRow{lightgray!40} 8443 & 7 & $L_{2}(8443^{2})$, $S_{4}(8443)$\\
8447 & 16 & $U_{3}(8447)$\\
\SetRow{lightgray!40} 8461 & 10 & $U_{3}(1777)$, $U_{4}(1777)$, $U_{3}(8461)$\\
8501 & 17 & $L_{4}(4481)$, $L_{2}(4481^{2})$, $S_{4}(4481)$, $U_{3}(8501)$\\
\SetRow{lightgray!40} 8521 & 9 & $L_{2}(8521^{2})$, $S_{4}(8521)$\\
8527 & 12 & $U_{3}(8527)$\\
\SetRow{lightgray!40} 8537 & 7 & $L_{2}(8537^{2})$, $S_{4}(8537)$, $U_{3}(8537)$, $U_{4}(8537)$\\
8539 & 6 & $U_{3}(5987)$\\
\SetRow{lightgray!40} 8563 & 12 & $L_{3}(8563)$\\
8573 & 13 & $L_{3}(8573)$, $L_{4}(8573)$, $L_{2}(8573^{2})$, $S_{4}(8573)$\\
\SetRow{lightgray!40} 8581 & 22 & $L_{2}(131^{2})$, $S_{4}(131)$, $U_{4}(131)$, $L_{2}(8581^{2})$, $S_{4}(8581)$\\
8597 & 5 & $L_{2}(8597^{2})$, $S_{4}(8597)$\\
\SetRow{lightgray!40} 8599 & 12 & $L_{3}(7393)$\\
8609 & 18 & $L_{2}(6779^{2})$, $S_{4}(6779)$, $U_{4}(6779)$\\
\SetRow{lightgray!40} 8623 & 6 & $U_{3}(8623)$\\
8629 & 17 & $U_{3}(3307)$, $U_{4}(3307)$, $U_{3}(5323)$, $U_{4}(5323)$\\
\SetRow{lightgray!40} 8641 & 11 & $L_{2}(1583^{2})$, $S_{4}(1583)$, $U_{4}(1583)$, $U_{3}(8641)$\\
8647 & 19 & $U_{3}(7853)$, $L_{3}(8647)$\\
\SetRow{lightgray!40} 8663 & 9 & $L_{2}(8663^{2})$, $S_{4}(8663)$\\
8669 & 16 & $L_{4}(4793)$, $L_{2}(4793^{2})$, $S_{4}(4793)$, $L_{2}(8669^{2})$, $S_{4}(8669)$, $U_{3}(8669)$, $U_{4}(8669)$\\
\SetRow{lightgray!40} 8681 & 11 & $L_{2}(3911^{2})$, $S_{4}(3911)$\\
8689 & 7 & $L_{2}(8689^{2})$, $S_{4}(8689)$\\
\SetRow{lightgray!40} 8707 & 11 & $L_{2}(8707^{2})$, $S_{4}(8707)$, $U_{3}(8707)$, $U_{4}(8707)$\\
8713 & 8 & $L_{3}(8713)$\\
\SetRow{lightgray!40} 8719 & 14 & $L_{3}(2281)$\\
8731 & 9 & $U_{3}(3659)$, $U_{4}(3659)$\\
\SetRow{lightgray!40} 8737 & 11 & $L_{2}(2269^{3})$, $G_2(2269)$, $U_{3}(2269)$, $U_{3}(6469)$, $U_{4}(6469)$, $U_{3}(8737)$\\
8747 & 8 & $L_{3}(8747)$\\
\SetRow{lightgray!40} 8761 & 24 & $L_{3}(1733)$, $L_{4}(1733)$, $L_{3}(7027)$, $L_{2}(8293^{2})$, $S_{4}(8293)$\\
8779 & 9 & $L_{3}(8779)$, $L_{2}(8779^{3})$, $G_2(8779)$, $U_{3}(8779)$\\
\SetRow{lightgray!40} 8783 & 23 & $L_{2}(8783^{2})$, $S_{4}(8783)$\\
8807 & 15 & $L_{2}(8807^{2})$, $S_{4}(8807)$\\
\SetRow{lightgray!40} 8821 & 15 & $U_{3}(467^{2})$, $U_{3}(2437)$, $U_{4}(2437)$, $U_{3}(8821)$\\
8831 & 8 & $U_{3}(8831)$\\
\SetRow{lightgray!40} 8837 & 4 & $U_{3}(8837)$\\
8839 & 13 & $U_{3}(4373)$, $U_{4}(4373)$\\
\SetRow{lightgray!40} 8861 & 7 & $L_{2}(8861^{2})$, $S_{4}(8861)$, $U_{3}(8861)$, $U_{4}(8861)$\\
8863 & 6 & $U_{3}(8863)$\\
\SetRow{lightgray!40} 8867 & 22 & $U_{3}(8867)$\\
8887 & 9 & $L_{2}(8887^{2})$, $S_{4}(8887)$\\
\SetRow{lightgray!40} 8893 & 35 & $L_{4}(2851)$, $L_{2}(2851^{2})$, $S_{4}(2851)$, $U_{3}(8893)$\\
8923 & 8 & $L_{3}(3847)$\\
\SetRow{lightgray!40} 8929 & 6 & $L_{3}(4339)$\\
8933 & 11 & $L_{2}(8171^{2})$, $S_{4}(8171)$\\
\SetRow{lightgray!40} 8941 & 15 & $L_{4}(5861)$, $L_{2}(5861^{2})$, $S_{4}(5861)$, $L_{3}(8941)$\\
8951 & 17 & $L_{3}(8951)$, $L_{2}(8951^{3})$, $G_2(8951)$, $U_{3}(8951)$\\
\SetRow{lightgray!40} 8971 & 31 & $L_{3}(8629)$, $U_{3}(8971)$\\
8999 & 7 & $L_{2}(8999^{2})$, $S_{4}(8999)$, $U_{3}(8999)$, $U_{4}(8999)$\\
\SetRow{lightgray!40} 9001 & 11 & $L_{4}(1237)$, $L_{2}(1237^{2})$, $S_{4}(1237)$, $L_{3}(9001)$\\
9007 & 6 & $L_{3}(9007)$\\
\SetRow{lightgray!40} 9011 & 5 & $L_{2}(9011^{2})$, $S_{4}(9011)$\\
9043 & 9 & $L_{2}(9043^{2})$, $S_{4}(9043)$\\
\SetRow{lightgray!40} 9049 & 14 & $L_{4}(7687)$, $L_{2}(7687^{2})$, $S_{4}(7687)$\\
9067 & 29 & $L_{3}(9067)$, $L_{4}(9067)$, $L_{2}(9067^{2})$, $S_{4}(9067)$\\
\SetRow{lightgray!40} 9091 & 22 & $L_{3}(3389^{2})$, $L_{2}(3389^{3})$, $S_{6}(3389)$, $O_{7}(3389)$, $O^+_{8}(3389)$, $G_2(3389)$, $U_{3}(3389)$, $U_{4}(3389)$, $L_{3}(9091)$\\
9103 & 14 & $L_{3}(4723)$, $L_{2}(4723^{3})$, $G_2(4723)$, $L_{3}(9103)$, $L_{2}(9103^{3})$, $G_2(9103)$, $U_{3}(9103)$\\
\SetRow{lightgray!40} 9109 & 24 & $U_{3}(3121)$, $L_{3}(9109)$, $L_{2}(9109^{3})$, $G_2(9109)$, $U_{3}(9109)$\\
9127 & 9 & $U_{3}(3011)$, $U_{4}(3011)$\\
\SetRow{lightgray!40} 9133 & 8 & $L_{3}(3797)$, $L_{2}(9133^{2})$, $S_{4}(9133)$\\
9157 & 8 & $L_{2}(2203^{2})$, $S_{4}(2203)$, $U_{4}(2203)$\\
\SetRow{lightgray!40} 9161 & 23 & $L_{4}(5^{5})$, $L_{2}(5^{10})$, $L_{3}(5^{10})$, $S_{4}(5^{5})$, $S_{6}(5^{5})$, $O_{7}(5^{5})$, $O^+_{8}(5^{5})$, $U_{5}(25)$, $U_{6}(25)$, $U_{4}(5^{5})$\\
9173 & 12 & $L_{4}(6659)$, $L_{2}(6659^{2})$, $S_{4}(6659)$\\
\SetRow{lightgray!40} 9181 & 13 & $L_{3}(1009)$, $L_{2}(1009^{3})$, $G_2(1009)$, $L_{3}(8171)$, $L_{4}(8171)$, $L_{3}(9181)$\\
9187 & 14 & $U_{3}(9187)$\\
\SetRow{lightgray!40} 9199 & 9 & $L_{3}(3767)$, $L_{3}(5431)$, $L_{2}(5431^{3})$, $G_2(5431)$\\
9209 & 18 & $L_{2}(8863^{2})$, $S_{4}(8863)$, $U_{4}(8863)$, $L_{2}(9209^{2})$, $S_{4}(9209)$\\
\SetRow{lightgray!40} 9221 & 9 & $L_{2}(9221^{2})$, $S_{4}(9221)$\\
9227 & 14 & $L_{3}(9227)$\\
\SetRow{lightgray!40} 9239 & 4 & $U_{3}(9239)$\\
9241 & 20 & $U_{3}(167)$, $U_{4}(167)$, $L_{3}(9241)$\\
\SetRow{lightgray!40} 9257 & 24 & $L_{2}(1097^{2})$, $S_{4}(1097)$, $U_{4}(1097)$\\
9277 & 11 & $L_{3}(601)$, $L_{4}(601)$, $L_{4}(8389)$, $L_{2}(8389^{2})$, $S_{4}(8389)$, $L_{3}(9277)$\\
\SetRow{lightgray!40} 9281 & 4 & $L_{3}(9281)$\\
9283 & 19 & $L_{3}(2843)$, $L_{4}(2843)$, $L_{3}(2843^{2})$, $L_{2}(2843^{3})$, $S_{6}(2843)$, $O_{7}(2843)$, $O^+_{8}(2843)$, $G_2(2843)$\\
\SetRow{lightgray!40} 9337 & 7 & $L_{3}(4937)$, $L_{4}(4937)$\\
9341 & 6 & $L_{4}(6703)$, $L_{2}(6703^{2})$, $S_{4}(6703)$\\
\SetRow{lightgray!40} 9343 & 10 & $L_{2}(6113^{3})$, $G_2(6113)$, $U_{3}(6113)$\\
9349 & 29 & $L_{2}(73^{6})$, $S_{4}(73^{3})$, $G_2(73^{2})$, $^3D_4(73)$, $U_{3}(73^{2})$, $U_{3}(4021)$\\
\SetRow{lightgray!40} 9377 & 18 & $L_{2}(6529^{2})$, $S_{4}(6529)$, $U_{4}(6529)$\\
9391 & 15 & $L_{3}(983)$, $L_{4}(983)$, $L_{3}(983^{2})$, $L_{2}(983^{3})$, $S_{6}(983)$, $O_{7}(983)$, $O^+_{8}(983)$, $G_2(983)$\\
\SetRow{lightgray!40} 9403 & 15 & $L_{3}(9403)$, $L_{2}(9403^{3})$, $G_2(9403)$, $U_{3}(9403)$\\
9419 & 4 & $L_{3}(9419)$\\
\SetRow{lightgray!40} 9421 & 12 & $L_{3}(9421)$\\
9431 & 5 & $L_{2}(9431^{2})$, $S_{4}(9431)$\\
\SetRow{lightgray!40} 9433 & 7 & $L_{2}(8419^{2})$, $S_{4}(8419)$\\
9439 & 25 & $L_{3}(733)$, $L_{4}(733)$\\
\SetRow{lightgray!40} 9461 & 6 & $L_{2}(7951^{2})$, $S_{4}(7951)$, $U_{3}(9461)$\\
9463 & 15 & $L_{3}(607)$, $L_{4}(607)$, $L_{3}(607^{2})$, $L_{2}(607^{3})$, $S_{6}(607)$, $O_{7}(607)$, $O^+_{8}(607)$, $G_2(607)$, $L_{2}(9463^{2})$, $S_{4}(9463)$\\
\SetRow{lightgray!40} 9467 & 9 & $L_{2}(9467^{2})$, $S_{4}(9467)$\\
9479 & 17 & $L_{2}(9479^{2})$, $S_{4}(9479)$, $U_{3}(9479)$, $U_{4}(9479)$\\
\SetRow{lightgray!40} 9491 & 8 & $L_{3}(9491)$\\
9511 & 12 & $U_{3}(3491)$\\
\SetRow{lightgray!40} 9551 & 38 & $U_{3}(9551)$\\
9587 & 19 & $L_{2}(9587^{2})$, $S_{4}(9587)$, $U_{3}(9587)$, $U_{4}(9587)$\\
\SetRow{lightgray!40} 9601 & 14 & $L_{3}(9601)$\\
9619 & 6 & $U_{3}(9619)$\\
\SetRow{lightgray!40} 9631 & 17 & $L_{3}(1621)$, $L_{4}(1621)$, $L_{3}(8009)$, $U_{3}(9631)$\\
9643 & 8 & $U_{3}(4597)$\\
\SetRow{lightgray!40} 9649 & 17 & $L_{3}(9649)$, $L_{2}(9649^{3})$, $G_2(9649)$, $U_{3}(9649)$\\
9661 & 20 & $L_{4}(139)$, $L_{2}(139^{2})$, $S_{4}(139)$\\
\SetRow{lightgray!40} 9677 & 4 & $L_{3}(9677)$\\
9679 & 12 & $U_{3}(9679)$\\
\SetRow{lightgray!40} 9689 & 11 & $L_{2}(7477^{2})$, $S_{4}(7477)$\\
9697 & 24 & $L_{3}(9697)$\\
\SetRow{lightgray!40} 9721 & 17 & $L_{3}(9721)$, $L_{2}(9721^{3})$, $G_2(9721)$, $U_{3}(9721)$\\
9733 & 8 & $U_{3}(9733)$\\
\SetRow{lightgray!40} 9739 & 12 & $L_{2}(6971^{3})$, $G_2(6971)$, $U_{3}(6971)$, $L_{3}(9739)$, $L_{4}(9739)$, $L_{2}(9739^{2})$, $S_{4}(9739)$\\
9791 & 14 & $U_{3}(9791)$\\
\SetRow{lightgray!40} 9811 & 13 & $U_{3}(9209)$, $U_{4}(9209)$, $L_{3}(9811)$, $L_{4}(9811)$, $L_{2}(9811^{2})$, $S_{4}(9811)$\\
9817 & 26 & $L_{2}(4027^{6})$, $S_{4}(4027^{3})$, $G_2(4027^{2})$, $^3D_4(4027)$, $U_{3}(4027^{2})$, $L_{3}(8861)$, $L_{4}(8861)$, $L_{3}(8861^{2})$, $L_{2}(8861^{3})$, $S_{6}(8861)$, $O_{7}(8861)$, $O^+_{8}(8861)$, $G_2(8861)$\\
\SetRow{lightgray!40} 9833 & 11 & $L_{2}(9833^{2})$, $S_{4}(9833)$, $U_{3}(9833)$, $U_{4}(9833)$\\
9839 & 15 & $L_{2}(9839^{2})$, $S_{4}(9839)$\\
\SetRow{lightgray!40} 9851 & 8 & $L_{3}(9851)$\\
9859 & 15 & $U_{3}(4751)$, $U_{4}(4751)$\\
\SetRow{lightgray!40} 9871 & 15 & $L_{2}(9871^{2})$, $S_{4}(9871)$\\
9883 & 6 & $L_{3}(9883)$\\
\SetRow{lightgray!40} 9887 & 19 & $L_{2}(9887^{2})$, $S_{4}(9887)$, $U_{3}(9887)$, $U_{4}(9887)$\\
9901 & 9 & $L_{2}(9901^{2})$, $S_{4}(9901)$\\
\SetRow{lightgray!40} 9907 & 19 & $L_{3}(6571)$, $L_{4}(6571)$\\
9923 & 8 & $U_{3}(9923)$\\
\SetRow{lightgray!40} 9929 & 7 & $L_{3}(9929)$, $L_{2}(9929^{3})$, $G_2(9929)$, $U_{3}(9929)$\\
9931 & 14 & $L_{3}(4231)$, $L_{2}(4231^{3})$, $G_2(4231)$\\
\SetRow{lightgray!40} 9941 & 11 & $L_{2}(9941^{2})$, $S_{4}(9941)$\\
9949 & 25 & $L_{2}(2543^{2})$, $S_{4}(2543)$, $L_{2}(9949^{2})$, $S_{4}(9949)$, $U_{3}(9949)$, $U_{4}(9949)$\\
\SetRow{lightgray!40} 9967 & 9 & $L_{3}(457)$, $L_{4}(457)$\\
\hline
\end{longtblr}

\end{document}